\documentclass{amsart}
\usepackage{amsfonts, amssymb}
\usepackage{amsmath}
\usepackage{amsthm}

\numberwithin{equation}{section}

\newtheorem{theorem}{Theorem}[section]

\newtheorem{corollary}[theorem]{Corollary}

\begin{document}

\title{Existence of meromorphic solutions of first order difference equations}

\author{Risto Korhonen}
\address{Department of Physics and Mathematics, University of Eastern Finland, P.O. Box 111,
FI-80101 Joensuu, Finland}
\email{risto.korhonen@uef.fi}
\thanks{The first author would like to thank the partial support by the Academy of Finland grant 286877.}

\author{Yueyang Zhang}
\address{Department of Physics and Mathematics, University of Eastern Finland, P.O. Box 111,
FI-80101 Joensuu, Finland}
\email{yueyang.zhang@uef.fi}
\thanks{The second author was supported in part by NSFC (No.~11371225) and the China Scholarship Council (CSC)}


\date{\today}

\commby{}

\begin{abstract}
It is shown that if
    \begin{equation}\label{abstract_eq}
    f(z+1)^n=R(z,f),\tag{\dag}
    \end{equation}
where $R(z,f)$ is rational in $f$ with meromorphic coefficients and $\deg_f(R(z,f))=n$, has an admissible meromorphic solution, then either $f$ satisfies a difference linear or Riccati equation with meromorphic coefficients, or \eqref{abstract_eq} can be transformed into one in a list of ten equations with certain meromorphic or algebroid coefficients. In particular, if \eqref{abstract_eq}, where the assumption $\deg_f(R(z,f))=n$ has been discarded, has rational coefficients and a transcendental meromorphic solution $f$ of hyper-order $<1$, then either $f$ satisfies a difference linear or Riccati equation with rational coefficients, or \eqref{abstract_eq} can be transformed into one in a list of five equations which consists of four difference Fermat equations and one equation which is a special case of the symmetric QRT map. Solutions to all of these equations are presented in terms of Weierstrass or Jacobi elliptic functions, or in terms of meromorphic functions which are solutions to a difference Riccati equation. This provides a natural difference analogue of Steinmetz' generalization of Malmquist's theorem.
\end{abstract}
\maketitle



\section{Introduction}\label{intro} 
Global existence of large classes of meromorphic solutions is a rare property for a differential equation to have. According to a classical result due to Malmquist~\cite{malmquist1913fonctions}, if the first order differential equation
    \begin{equation}\label{first_order_de}
    f'=R(z,f),
    \end{equation}
where $R(z,f)$ is rational in both arguments, has a transcendental meromorphic solution, then \eqref{first_order_de} reduces into the Riccati equation
    \begin{equation}\label{riccati}
    f'=a_2 f^2 + a_1 f + a_0
    \end{equation}
with rational coefficients $a_0$, $a_1$ and $a_2$. Generalizations of Malmquist's theorem for the equation
    \begin{equation}\label{first_order_de_n}
    (f')^n=R(z,f),\qquad n\in\mathbb{N},
    \end{equation}
have been given by Yosida \cite{yosida:33} and Laine \cite{laine:71}. Steinmetz \cite{steinmetz:78}, and Bank and Kaufman \cite{bank1980growth} proved that if \eqref{first_order_de_n} has rational coefficients and a transcendental meromorphic solution, then by a suitable M\"obius transformation, \eqref{first_order_de_n} can be either mapped to \eqref{riccati}, or to one of the equations in the following list:
    \begin{equation*}
    \begin{split}
    (f')^2 &= a(f-b)^2(f-\tau_1)(f-\tau_2)\\
    (f')^2 &= a(f-\tau_1)(f-\tau_2)(f-\tau_3)(f-\tau_4)\\
    (f')^3 &= a(f-\tau_1)^2(f-\tau_2)^2(f-\tau_3)^2\\
    (f')^4 &= a(f-\tau_1)^2(f-\tau_2)^3(f-\tau_3)^3\\
    (f')^6 &= a(f-\tau_1)^3(f-\tau_2)^4(f-\tau_3)^5\\
    \end{split}
    \end{equation*}
where $a$ and $b$ are rational functions, and $\tau_1,\ldots,\tau_4$ are distinct constants.

The existence of globally meromorphic solutions is somewhat more common in the case of difference equations, as compared to differential equations. It was shown by Shimomura \cite{shimomura:81} that the difference equation
    \begin{equation*}
    f(z+1) = P(f(z)),
    \end{equation*}
where $P(f(z))$ is a polynomial in $f(z)$ with constant coefficients, always has a non-trivial entire solution. On the other hand, Yanagihara \cite{yanagihara:80} showed that the difference equation
    \begin{equation*}
    f(z+1) = R(f(z)),
    \end{equation*}
where $R(f(z))$ is rational in $f(z)$ having constant coefficients, has a non-trivial meromorphic solution no matter how $R$ is chosen. Yanagihara \cite{yanagihara:85} also considered higher order equations and showed, for instance, that the difference equation
    \begin{equation*}
    \alpha_{n}f(z+n)+\alpha_{n-1}f(z+n-1)+\cdots+\alpha_{1}f(z+1)=R(f(z)), \qquad \alpha_1,\ldots,\alpha_n\in\mathbb{C},
    \end{equation*}
has a non-trivial meromorphic solution if the degree $p$ of the numerator $P(f(z))$ of the rational function $R(f(z))$ satisfies $p\geq q+2$,  where $q$ is the degree of the denominator $Q(f(z))$, and $P(f(z))$ and $Q(f(z))$ have no common factors.

Ablowitz, Halburd and Herbst \cite{AblowitzHalburdHerbst2000} suggested that the existence of sufficiently many finite-order meromorphic solutions of a difference equation is a good difference analogue of the Painlev\'e property for differential equations. An ordinary differential equation is said to have the Painlev\'e property when all solutions are single-valued around all movable singularities. They showed, for instance, that if the difference equation
    \begin{equation}\label{ahh_eq}
    f(z+1)+f(z-1) = R(z,f(z)),
    \end{equation}
where $R(z,f(z))$ is rational in both arguments, has a transcendental meromorphic solution of finite order, then $\deg_f(R(z,f(z)))\leq 2$. Their results are consistent with Yanagihara's work on the first order equation \cite{yanagihara:80} in which he proved that if
    \begin{equation}\label{yanagihara_eq}
    f(z+1)= R(z,f(z)),
    \end{equation}
where $R(z,f(z))$ is rational in both arguments, has a transcendental meromorphic solution of hyper-order strictly less than one, then $\deg_f(R(z,f(z)))= 1$ and thus \eqref{yanagihara_eq} reduces into the difference Riccati equation. This is a natural difference analogue of Malmquist's 1913 result on differential equations.  Halburd and the first author \cite{halburdrk:LMS2006} showed that if \eqref{ahh_eq}, where the right hand side has meromorphic coefficients, has an admissible meromorphic solution $f$ of finite order, then either $f$ satisfies a difference Riccati equation, or a linear transformation of \eqref{ahh_eq} reduces it into one in a short list of difference equations which consists solely of difference Painlev\'e equations and equations related to them, linear equations and linearizable equations.
The finite-order condition was relaxed into hyper-order strictly less than one by Halburd, the first author and Tohge \cite{halburdkt:14TAMS}.

The purpose of this paper is to present a natural difference analogue of Steinmetz' generalization of Malmquist's theorem. We will show that if the difference equation
    \begin{equation}\label{intro_eq}
    f(z+1)^n=R(z,f),
    \end{equation}
with rational coefficients has a transcendental meromorphic solution $f$ of hyper-order $<1$, then either $f$ satisfies a difference linear or Riccati equation
    \begin{eqnarray}
    f(z+1)&=& a_1(z)f(z)+a_2(z),\label{lineareq_intro}\\
    f(z+1)&=&\frac{b_1(z)f(z)+b_2(z)(z)}{f(z)+b_3(z)}\label{driccati0_intro},
     \end{eqnarray}
where $a_i(z)$ and $b_j(z)$ are rational functions, or, by implementing a transformation $f\rightarrow \alpha f$ or $f\rightarrow 1/(\alpha f)$ with an algebraic function $\alpha$ of degree at most~2, \eqref{intro_eq} reduces into one of the following equations:
   \begin{eqnarray}
   f(z+1)^2 &=& 1-f(z)^2,\label{yanagiharaeq11_intro}\\
   f(z+1)^2 &=& 1-\left(\frac{\delta(z) f(z)-1}{f(z)-\delta(z)}\right)^2,\label{yanagiharaeq12_intro}\\
   f(z+1)^{2} &=& 1-\left(\frac{f(z)+3}{f(z)-1}\right)^{2}\label{yanagiharaeq13_intro}\\
   f(z+1)^{2} &=& \frac{f(z)^2 - \kappa^2}{f(z)^2-1}, \label{yanagiharaeq14_intro}\\
   f(z+1)^{3} &=& 1-f(z)^{-3},\label{yanagiharaeq15_intro}
  \end{eqnarray}
where $\delta(z)$ $(\not\equiv\pm1)$ is an algebraic function of degree $2$ at most and $\kappa^2\not=0,1$ is a constant. Equations \eqref{lineareq_intro} and \eqref{driccati0_intro} are the linear and the difference Riccati equation, respectively, equations \eqref{yanagiharaeq11_intro}--\eqref{yanagiharaeq13_intro} and \eqref{yanagiharaeq15_intro} are difference Fermat equations, while \eqref{yanagiharaeq14_intro} is a special case of the symmetric QRT map. We will present finite-order meromorphic solutions to autonomous versions of all of these equations in terms of Weierstrass or Jacobi elliptic functions, or in terms of meromorphic functions which are solutions to a difference Riccati equation, in Section~\ref{malm_sec} below.

Under the condition that the meromorphic solution $f$ of \eqref{intro_eq} is of hyper-order less than~1, it can actually be shown that $\deg_f(R(z,f))=n$ by using an asymptotic relation between the Nevanlinna characteristics (see Section~\ref{prelim_sec} below) $T(r,f(z+1))$ and $T(r,f(z))$ from \cite{halburdkt:14TAMS}, and an identity due to Valiron \cite{Valiron1931} (see also \cite{Laine1993}). By discarding the assumption that the meromorphic solution is of hyper-order $<1$, and considering the more general case of admissible meromorphic solutions of \eqref{intro_eq} with meromorphic coefficients and such that $\deg_f(R(z,f))=n$, it follows either that $f$ satisfies \eqref{lineareq_intro} or \eqref{driccati0_intro}, or \eqref{intro_eq} can be transformed into one of the equations \eqref{yanagiharaeq11_intro}--\eqref{yanagiharaeq15_intro}, but now with meromorphic coefficients $a_i(z)$ and $b_j(z)$, an algebroid function $\delta(z)=\delta_2(z)$ of degree at most~2 and with $\kappa^2=\kappa_1(z)^2$ being a meromorphic periodic function of period $1$, or \eqref{intro_eq} becomes one of the following equations:
   \begin{eqnarray}
   f(z+1)^2 &=& \delta_1(z)(f(z)^2-1), \label{intro_eq_list2_1}\\
   f(z+1)^2 &=& \delta_3(z) (1-f(z)^{-2}),\\
   f(z+1)^2 &=& \frac{\kappa_2(z+1)^2f(z)^2-1}{f(z)^2-1},\\
   f(z+1)^2 &=& \theta\frac{f(z)^2-\kappa_3(z)f(z)+1}{f(z)^2+\kappa_3(z)f(z)+1},\\
   f(z+1)^3 &=& 1-f(z)^3,\label{intro_eq_list2_5}
  \end{eqnarray}
where $\theta=\pm1$ and $\delta_1(z)$, $\delta_3(z)$, $\kappa_2(z)^2$, $\kappa_3(z)^2$ are meromorpic functions each of which satisfies a certain difference equation, see Section~\ref{malm_sec} below. In particular, if the coefficients of \eqref{intro_eq} are rational functions, then $\theta=\pm1$ and $\delta_1(z)$, $\delta_3(z)$, $\kappa_2(z)^2$, $\kappa_3(z)^2$ are all constants. We will show that meromorphic solutions of autonomous versions of \eqref{intro_eq_list2_1}--\eqref{intro_eq_list2_5} can be characterized by Weierstrass or Jacobi elliptic functions composed with certain entire functions, but none of them is of hyper-order $<1$.

The remainder of this paper is organized in the following way. Section~\ref{prelim_sec} contains the necessary notations and preliminary results needed to prove our main theorem given in Section~\ref{malm_sec}. Section~\ref{malm_sec} also contains a discussion of the equations \eqref{yanagiharaeq11_intro}--\eqref{intro_eq_list2_5}. In particular, we will show that either their autonomous versions have meromorphic solutions expressible in terms of Weierstrass or Jacobi elliptic functions, or solutions of these equations can be explicitly expressed in terms of meromorphic functions which are solutions of certain difference Riccati equations. The proof of our main result, Theorem~\ref{dtheorem} below, has been split into the three remaining sections \ref{proof1_sec}--\ref{proof3_sec}.

\section{Preliminaries}\label{prelim_sec}

In this section, we introduce some preliminaries for the proof of our main results. In the following, a meromorphic function is always meromorphic on the whole complex plane $\mathbb{C}$. Let $\mathcal{M}$ denote the field of meromorphic functions and let $f(z)\in\mathcal{M}$.
We assume that the readers are familiar with the basic notations and fundamental results of Nevanlinna theory \cite{Hayman1964Meromorphic}, such as $m(r,f)$, $N(r,f)$, $T(r,f)$ and Nevanlinna's main theorems.
Moreover, we use the notation $\sigma(f)$ to denote the \emph{order of growth} of $f(z)$ which is defined to be
\begin{equation*}
\sigma(f)=\limsup_{r\to\infty}\frac{\log T(r,f)}{\log{r}},
\end{equation*}
and, if $f(z)$ is of order $\sigma(f)=\infty$, we use the notation $\sigma_2(f)$ to denote the \emph{hyper-order} of $f(z)$ which is defined to be
\begin{equation*}
\sigma_2(f)=\limsup_{r\to\infty}\frac{\log\log T(r,f)}{\log{r}}.
\end{equation*}
Let $S(r,f)$ denote any quantity that satisfies the condition $S(r,f)=o(1)T(r,f)$, $r\rightarrow\infty$, outside of a possible exceptional set of finite logarithmic measure.
For a meromorphic function $c(z)\in\mathcal{M}$, if $c(z)$ satisfies $T(r,c(z))=S(r,f)$, then $c(z)$ is said to be \emph{small} compared to $f(z)$. For example, all rational functions are small with respect to any transcendental meromorphic function.
Denote by $\mathcal{S}(f)$ the field of all small functions of $f(z)$, i.e.,
\begin{equation*}
\mathcal{S}(f)=\{c(z)\in\mathcal{M}: T(r,c(z))=S(r,f)\}.
\end{equation*}
Set $\hat{\mathcal{S}}(f)=\mathcal{S}(f)\cup\{\infty\}$. A meromorphic solution $f(z)$ of a differential equation is called \emph{admissible} if all coefficients of the equation are in $\mathcal{S}(f)$ \cite{Laine1993}. This definition is very natural for differential equations, but for difference equations to be studied in this paper we need to be slightly more careful due to the possible effect of the shift operation to the growth of the solution and of the coefficients. We define the following set for small functions of a meromorphic function $f(z)$
\begin{equation*}
\mathcal{S}'(f)=\{c(z)\in\mathcal{M}: T(r,c(z+n))=S(r,f(z+n)), n\in\mathbb{N}\},
\end{equation*}
and in what follows we say that a meromorphic solution $f(z)$ of a difference equation is admissible if all coefficients of the equation are in $\mathcal{S}'(f)$.
If the hyper-order of $f(z)$ is less than one, then by \cite[Lemma~8.3]{halburdkt:14TAMS} it follows that
     \begin{equation}\label{T shift}
    T(r,f(z+1)) = T(r,f) + S(r,f),
    \end{equation}
and so in this case shifting a difference equation does not affect the admissibility of solutions.

The fundamental results of Nevanlinna's theory are what are known as the First Main Theorem and the Second Main Theorem. Nevanlinna's Second Main Theorem can be simplified into the following form:
    \begin{equation}\label{secondNevanlinna}
    (q-2)T(r,f)\leq \sum_{i=1}^q\overline{N}(r,a_i,f)+S(r,f),
    \end{equation}
where $q\geq 3$, $a_i\in\mathbb{C}\cup\{\infty\}$, and $\overline{N}(r,a_i,f)$ is the truncated counting function for $a_i$-points (poles) of $f(z)$, and the error term here $S(r,f)=O(\log r)$, $r\rightarrow\infty$, when $f(z)$ is of finite order and $S(r,f)=O(\log rT(r,f))$, $r\rightarrow\infty$, outside a possible set of finite linear measure when $f(z)$ is of infinite order. Yamanoi \cite[Corollary~1]{yamanoi:04} generalized Nevanlinna's Second Main Theorem by proving that, for any collection of functions $c_i(z)\in\hat{\mathcal{S}}(f)$, $i=1,\ldots,q$, of $f(z)$, the following inequality
     \begin{equation}\label{Yamanoi co}
    (q-2-\varepsilon)T(r,f)\leq \sum_{i=1}^q\overline{N}(r,c_i,f), \quad \text{for all} \quad \varepsilon>0,
    \end{equation}
holds outside an exceptional set of $E\subset (0,\infty)$ satisfying $\int_{E}d\log\log r<\infty$. The inequality \eqref{Yamanoi co} implies a generalization of the defect relation for small functions of $f(z)$ \cite[Corollary~1]{yamanoi:04} and that
     \begin{equation}\label{de relation}
     \sum_{i=1}^q\Theta(c_i,f)\leq 2,
     \end{equation}
where the quantity $\Theta(c,f)$, $c=c(z)\in\hat{\mathcal{S}}(f)$, is defined to be
     \begin{equation*}
     \Theta(c,f)=1-\limsup_{r\to\infty}\frac{\overline{N}(r,c,f)}{T(r,f)}.
     \end{equation*}
Recall that a value $a$ is said to be a \emph{completely ramified value} of $f(z)$ when $f(z)-a=0$ has no simple roots. In what follows, we say that $c(z)\in\hat{\mathcal{S}}(f)$ is a completely ramified small function of $f(z)$ when $f(z)-c(z)=0$ has at most $S(r,f)$ many simple roots and that $c(z)\in\hat{\mathcal{S}}(f)$ is a Picard exceptional small function of $f(z)$ when $N(r,c,f)=S(r,f)$. A nonconstant meromorphic function $f(z)$ can have at most two Picard exceptional small functions. Moreover, by \eqref{de relation}, we have the following result.

\begin{theorem}\label{completelyrm}
A nonconstant meromorphic function $f(z)$ can have at most four completely ramified small functions.
\end{theorem}

When considering the meromorphic solution $f$ of \eqref{intro_eq}, we will do a transformation to $f$ using some algebroid functions and end up in a situation such that the considered function have some finite-sheeted branching. The classical version of Nevanlinna theory introduced above cannot be used to handle this situation. In this case, we need the Selberg-Ullrich theory, the algebroid version of Nevanlinna theory (see, for instance, \cite{Katajamaki1993algebroid}), which studies meromorphic functions on a finitely sheeted Riemann surface. All algebroid functions we need to consider in this paper are either small functions with respect to an admissible meromorphic solution $f$ of \eqref{intro_eq}, or can be obtained from it by a M\"obius transformation with small algebroid coefficients. Such functions could be described as ``almost meromorphic'' in the sense of Nevanlinna theory, since the presence of branch points actually only affects the small error term $S(r,f)$ in any of the estimates involving Nevanlinna functions. Correspondingly, $T(r,f)$ and $N(r,f)$ will denote the characteristic and counting functions of a finite-sheeted algebroid function $f$, and similarly with the rest of the Nevanlinna functions involving $f$.

For simplicity, from now on we will use the suppressed notations: $f=f(z)$, $\overline{f}=f(z+1)$ and $\underline{f}=f(z-1)$ for a meromorphic, or algebroid, function $f(z)$.

\section{Extension of the difference Malmquist theorem}\label{malm_sec}

As mentioned in the introduction, from equation \eqref{yanagihara_eq} Yanagihara \cite{yanagihara:80} obtained a difference analogue of Malmquist's theorem. Let $n\in\mathbb{N}$.  Nakamura and Yanagihara \cite{NakamuraYanagihara1989difference}, and Yanagihara \cite{Yanagihara1989difference}, considered meromorphic solutions of a more general case
    \begin{equation}\label{yanagiharaeq n}
   \overline{f}^n=P(z,f),
    \end{equation}
where $P(z,f)$ is a polynomial in $f$ with constants as coefficients and, in particular, Yanagihara \cite{Yanagihara1989difference} showed that if \eqref{yanagiharaeq n} has a meromorphic solution of finite order, then either $\overline{f}=Af+B$ for some constants $A,B$ or \eqref{yanagiharaeq n} reduces into $\overline{f}^2=1-f^2$. We consider the difference equation
    \begin{equation}\label{yanagiharaeq2}
   \overline{f}^n=R(z,f),
    \end{equation}
where $R(z,f)$ is rational in $f$ having meromorphic coefficients. If equation \eqref{yanagiharaeq2} has an admissible meromorphic solution $f$ such that the hyper-order $\sigma_2(f)$ of $f$ satisfies $\sigma_2(f)<1$, then, by \eqref{yanagiharaeq2}, it follows from \eqref{T shift} and the Valiron's identity \cite{Valiron1931} (see also \cite{Laine1993}) that
    \begin{equation*}
    \begin{split}
    \deg_f(R(z,f))T(r,f) &= T(r,R(z,f)) + S(r,f)\\
    &= T(r,\overline{f}^n) + S(r,f) \\
    &= nT(r,\overline{f}) + S(r,f) \\
    &= nT(r,f) + S(r,f),
    \end{split}
    \end{equation*}
which implies that $\deg_f(R(z,f))=n$. Assume that $\deg_f(R(z,f))=n$ and that all the coefficients of \eqref{yanagiharaeq2} are in $\mathcal{S}'(f)$. Then the characteristic function of $f(z)$ satisfies $T(r,f(z+j))=T(r,f(z))+S(r,f)$, $j\in\mathbb{N}$. Under these conditions, but with discarding the assumption $\sigma_2(f)<1$, we now use Nevanlinna theory to simplify equation \eqref{yanagiharaeq2}. We state our main results as follows, which are further extensions of the results on \eqref{yanagiharaeq n} of Yanagihara \cite{Yanagihara1989difference}.

\begin{theorem}\label{dtheorem}
Let $n\in\mathbb{N}$. If the difference equation \eqref{yanagiharaeq2} with $\deg_f(R(z,f))=n$ has an admissible meromorphic solution, then either $f$ satisfies a difference linear or Riccati equation:
    \begin{eqnarray}
    \overline{f}&=& a_1f+a_2,\label{lineareq}\\
    \overline{f}&=&\frac{b_1f+b_2}{f+b_3}\label{driccati0},
     \end{eqnarray}
where $a_i,b_j$ are small meromorphic functions; or, by a transformation $f\rightarrow \alpha f$ or $f\rightarrow 1/(\alpha f)$ with a small algebroid function $\alpha$ of degree at most~3, \eqref{yanagiharaeq2} reduces into one of the following equations:
   \begin{eqnarray}
   \overline{f}^2 &=& 1-f^2,\label{yanagiharaeq11}\\
   \overline{f}^2 &=& \delta_1(f^2-1),\label{yanagiharaeq11a}\\
   \overline{f}^2 &=& 1-\left(\frac{\delta_2 f-1}{f-\delta_2}\right)^2,\label{yanagiharaeq12}\\
   \overline{f}^2 &=& \delta_3(1-f^{-2}),\label{yanagiharaeq12a}\\
   \overline{f}^{2} &=& 1-\left(\frac{f+3}{f-1}\right)^{2},\label{yanagiharaeq13}\\
   \overline{f}^2 &=& \frac{f^2-\kappa_1^2}{f^2-1},\label{yanagiharaeq13a}\\
   \overline{f}^2 &=& \frac{\overline{\kappa}_2^2f^2-1}{f^2-1},\label{yanagiharaeq13b}\\
   \overline{f}^2 &=& \theta\frac{f^2-\kappa_3f+1}{f^2+\kappa_3f+1},\label{yanagiharaeq13c}\\
   \overline{f}^{3} &=& 1-f^3,\label{yanagiharaeq14}\\
   \overline{f}^{3} &=& 1-f^{-3},\label{yanagiharaeq15}
  \end{eqnarray}
where $\theta=\pm1$, $\delta_2\not\equiv\pm1$ is a small algebroid function of degree at most~2 and $\delta_1$, $\delta_3$, $\kappa_1^2$, $\kappa_2^2$, $\kappa_3^2$ are all small meromorphic functions satisfying $\overline{\delta}_1(\delta_1+1)+1=0$, $\overline{\delta}_3\delta_3=\overline{\delta}_3+\delta_3$, $\overline{\kappa}_1^2=\kappa_1^2$, $\overline{\kappa}_2^2\kappa_2^2=1$ and $\overline{\kappa}_3^2(\kappa_3^2-4)=2(1-\theta)\kappa_3^2-8(1+\theta)$.
\end{theorem}

Assume that in Theorem~\ref{dtheorem} all the coefficients of \eqref{yanagiharaeq2} are rational functions. Then it is seen that the coefficients $\delta_1,\delta_3,\kappa_1^2,\kappa_2^2,\kappa_3^2$ are all constants. In this case, it is easy to show by iteration and substitution that solutions to \eqref{yanagiharaeq11}, \eqref{yanagiharaeq11a}, \eqref{yanagiharaeq12a}, \eqref{yanagiharaeq13a}, \eqref{yanagiharaeq13b}, \eqref{yanagiharaeq14} and \eqref{yanagiharaeq15} are all periodic functions, which imply that for these equations the coefficient $\alpha$ in the transformation $f\rightarrow \alpha f$ or $f\rightarrow 1/(\alpha f)$ has no branch points. Also, from the proof of Theorem~\ref{dtheorem} we see that when obtaining \eqref{yanagiharaeq13} and \eqref{yanagiharaeq13c} the coefficient $\alpha$ in the transformation $f\rightarrow \alpha f$ or $f\rightarrow 1/(\alpha f)$ is rational. Thus the algebraic case of $\alpha$ can only occur when obtaining equation \eqref{yanagiharaeq12}.

In the rational coefficients case of \eqref{yanagiharaeq2}, solutions to \eqref{yanagiharaeq14} and \eqref{yanagiharaeq15} can be characterized by Weierstrass elliptic functions, solutions to equations \eqref{yanagiharaeq11a}, \eqref{yanagiharaeq12a}, \eqref{yanagiharaeq13a}--\eqref{yanagiharaeq13c} can be characterized by Jacobi elliptic functions and equations \eqref{yanagiharaeq11}, \eqref{yanagiharaeq12} and \eqref{yanagiharaeq13} can be explicitly solved in terms of functions which are solutions of certain difference Riccati equations, as is shown below.

Equations \eqref{yanagiharaeq11}, \eqref{yanagiharaeq12}, \eqref{yanagiharaeq13}, \eqref{yanagiharaeq14} and \eqref{yanagiharaeq15} are so-called \emph{Fermat} difference equations. In general, a Fermat equation is a function analogue of the Fermat diophantine equation $x^n+y^n=1$, i.e., $h(z)^n+g(z)^n=1$, where $n\geq2$ is an integer. Meromorphic solutions to Fermat equations have been clearly characterized, see \cite{Baker1966,Gross1966erratum,Gross1966}, for example. In particular, when $n=3$, all meromorphic solutions can be represented as: $h=H(\varphi)$, $g=\eta G(\varphi)=\eta H(-\varphi)=H(-\eta^2\varphi)$, where $\varphi=\varphi(z)$ is an entire function and $\eta$ is a cubic root of~1, and
    \begin{equation}\label{Weierstrassellptic0}
   H(z)=\frac{1+\wp'(z)/\sqrt{3}}{2\wp(z)}, \quad G(z)=\frac{1-\wp'(z)/\sqrt{3}}{2\wp(z)},
   \end{equation}
is a pair of solutions of the Fermat equation with $\wp(z)$ being the particular Weierstrass elliptic function
satisfying $\wp'(z)^2=4\wp(z)^3-1$. Now, for equation \eqref{yanagiharaeq14} the solution takes the form: $f=H(\varphi)$, $\overline{f}=\eta G(\varphi)=H(-\eta^2\varphi)$.
It follows that $\overline{\varphi}=-\eta^2\varphi+\omega$, where $\omega$ is a period of $\wp(z)$, which implies that $\varphi$ is a transcendental entire function of order at least~1. From \eqref{yanagiharaeq14} we have $T(r,\overline{f})=T(r,f)+O(1)$. Moreover, we have $\eta f+\overline{f}=\eta/\wp(\varphi)$, which yields $T(r,\wp(\varphi))\leq 2T(r,f)+O(1)$. By taking the derivatives on both sides of $\eta f+\overline{f}=\eta/\wp(\varphi)$ and combining it with the resulting equation, we get
    \begin{equation}\label{Weierstrassellptic011}
    \varphi'=-\frac{(\eta f+\overline{f})'}{\eta f+\overline{f}}\cdot\frac{\wp(\varphi)}{\wp'(\varphi)}.
    \end{equation}
For the elliptic function $\wp(z)$, we have
$N(r,\wp(\varphi))=T(r,\wp(\varphi))+O(\log rT(r,\wp(\varphi)))$ and it follows that $m(r,\wp(\varphi))=O(\log rT(r,\wp(\varphi)))=O(\log rT(r,f))$ and, similarly,
$m(r,1/\wp'(\varphi))=O(\log rT(r,f))$. Thus, by taking the proximity functions on both sides of \eqref{Weierstrassellptic011}, we get
    \begin{equation*}
    m(r,\varphi')\leq O(\log rT(r,f))+m(r,\wp(\varphi))+m\left(r,\frac{1}{\wp'(\varphi)}\right)+O(1)=O(\log rT(r,f))+O(1).
    \end{equation*}
Then, $T(r,\varphi')=m(r,\varphi')=O(\log rT(r,f))+O(1)$, which together with \cite[Lemma~1.1.2]{Laine1993} yields that $\sigma_2(f)\geq 1$ and so
equation \eqref{yanagiharaeq14} cannot have meromorphic solutions of hyper-order less than~1.

We now show that the autonomous versions of \eqref{yanagiharaeq11a}, \eqref{yanagiharaeq12a}, \eqref{yanagiharaeq13b} and \eqref{yanagiharaeq13c} cannot admit meromorphic solutions of hyper-order $<1$ either. We only need to consider equations \eqref{yanagiharaeq11a} and \eqref{yanagiharaeq13b} since \eqref{yanagiharaeq12a} and \eqref{yanagiharaeq13c} can be transformed into \eqref{yanagiharaeq11a} or \eqref{yanagiharaeq13b} in the following way. For equation \eqref{yanagiharaeq12a}, $\delta_3=2$ and we see that $f^2-1$ is written as $f^2-1=-(i\overline{f}f/\sqrt{2})^2=-f_0^2$. Then we have $T(r,f)=T(r,f_0)+O(1)$ and by substituting $f^2=-f_0^2+1$ into \eqref{yanagiharaeq12a} we get the autonomous case of equation \eqref{yanagiharaeq13b}. For equation \eqref{yanagiharaeq13c}, if we let $w=f+1/f$, then $w$ satisfies $T(r,w)=2T(r,f)+O(1)$ by Valiron's identity \cite{Valiron1931} (see also \cite{Laine1993}) and it follows that
   \begin{equation*}
   \overline{f}^2=\theta\frac{w-\kappa_3}{w+\kappa_3}, \quad \frac{1}{\overline{f}^2}=\theta\frac{w+\kappa_3}{w-\kappa_3}.
    \end{equation*}
The above two equations yield
    \begin{equation*}
    \overline{w}^2=\frac{2(\theta+1)w^2+2(\theta-1)\kappa_3^2}{w^2-\kappa_3^2}.
    \end{equation*}
If $\theta=1$, then we have $\overline{\kappa}_3^2(\kappa_3^2-4)=-16$ and by doing a transformation $w\rightarrow \kappa_3/w$ we get $\overline{w}^2=d_1(w^2-1)$,
which is the equation \eqref{yanagiharaeq11a} since $d_1=-\overline{\kappa}_3^2/4$ satisfies $\overline{d}_1(d_1+1)+1=0$; if $\theta=-1$, then we have $\overline{\kappa}_3^2(\kappa_3^2-4)=4\kappa_3^2$ and by doing a transformation $w\rightarrow \kappa_3/w$ we get $\overline{w}^2=d_2(1-w^{-2})$,
which is the equation \eqref{yanagiharaeq12a} since $d_2=\overline{\kappa}_3^2/4$ satisfies $\overline{d}_2d_2=\overline{d}_2+d_2$.

Recall that the Jacobi elliptic function $\text{sn}(z,k)$ with the elliptic modulus $k\in(0,1)$ satisfies the differential equation $\text{sn}'(z)^2=(1-\text{sn}(z)^2)(1-k^2\text{sn}(z)^2)$. For equation \eqref{yanagiharaeq11a}, $\delta_1$ is now a cubic root of~1 and we see that $f$ is twofold ramified over~$\pm1$. Denote $\delta_1=\eta^2$, where $\eta$ is the cubic root of unity and satisfies $\eta^2(\eta^2+1)+1=\eta^2+\eta+1=0$. Let $z_0$ be such that $f(z_0)=\pm i\eta$, by \eqref{yanagiharaeq11a} we have $f(z_0+1)=\pm1$, and so $f$ is also twofold ramified at $\pm i\eta$-points. Let $\varrho_1=\exp(7i\pi/12)$ and denote $\tau_1=(1-\varrho_1)/(1+\varrho_1)$. Then
    \begin{equation*}
    g_1=\tau_1\frac{f+\varrho_1}{f-\varrho_1}
    \end{equation*}
is twofold ramified over each of $\pm1,\pm\tau_1^2$, where $-\tau_1^2=\tan(7\pi/24)^2>1$.
Let $\text{sn}(\phi)=\text{sn}(\phi,-1/\tau_1^2)$ and $\phi_0$ be such that $\text{sn}(\phi_0)\not=\pm1,\pm\tau_1^2$. Letting $z_0$ be such that $g_1(z_0)=\text{sn}(\phi_0)$, it follows that there is a neighbourhood $U$ of $z_0$ such that $\phi_1(z)=\text{sn}^{-1}(g_1(z))$ is defined and holomorphic in $U$. By following the reasoning in the proof of \cite[Lemma~4.1]{Yanagihara1989difference}, we know that $\phi_1(z)$ can be continued analytically throughout the complex plane to an entire function. Therefore, $g_1$ is written as $g_1(z)=\text{sn}(\phi_1(z))$, where $\text{sn}(\phi_1)=\text{sn}(\phi_1,-1/\tau_1^2)$ satisfies the differential equation $\text{sn}'(\phi_1)^2=(1-\text{sn}(\phi_1)^2)(1-\text{sn}(\phi_1)^2/\tau_1^4)$. By taking the derivative of $g_1$ and combining it with the resulting equation, we get
    \begin{equation}\label{PPreq12aa1r3 co}
    \phi_1'^2=\frac{\text{sn}(\phi_1)^2}{\text{sn}'(\phi_1)^2}\cdot\left(\frac{g_1'}{g_1}\right)^2=\frac{g_1^2}{(1-g_1^2)(1-g_1^2/\tau_1^4)}\cdot\left(\frac{g_1'}{g_1}\right)^2.
    \end{equation}
Since $g_1$ is an elliptic function, we have
$m(r,g_1)=O(\log rT(r,f))$ and $m(r,1/(g_1-\varrho))=O(\log rT(r,f))$, $\varrho=\pm1,\pm\tau_1^2$. Note that $T(r,f)=T(r,g_1)+O(1)$. Now taking the proximity function on both sides of \eqref{PPreq12aa1r3 co} gives $T(r,\phi_1'^2)=m(r,\phi_1'^2)=O(\log rT(r,f))+O(1)$, which together with \cite[Lemma~1.1.2]{Laine1993} yields that $\phi_1'$ has order of growth strictly less than~1 when $\sigma_2(f)<1$. By iterating \eqref{yanagiharaeq11a} we obtain $f(z+3)^2=f(z)^2$ and so $f(z+6)=f(z)$ and it follows that $\text{sn}(\phi_1(z+6))=\text{sn}(\phi_1(z))$ giving $\phi_1(z+6)=\phi_1(z)+K_1$ with a period $K_1$, which is possible only when $\phi_1$ is a polynomial of degree~1 since $\phi_1'(z+6)=\phi_1'(z)$. This implies that $f$ is of finite order. But according to \cite[Lemma~9.1]{Yanagihara1989difference}, \eqref{yanagiharaeq11a} cannot admit any meromorphic solutions of finite order when $\delta_1$ is a constant. Therefore, \eqref{yanagiharaeq11a} cannot admit any meromorphic solutions of hyper-order strictly less than~1 when $\delta_1$ is a constant.

For equation \eqref{yanagiharaeq13b}, if $\kappa_2$ is a constant, then $\kappa_2^2=-1$ and we see that $f$ is twofold ramified over each of $\pm i,\pm1$. Let $\varrho_2=\exp(i\pi/4)$ and denote $\tau_2=(1+\varrho_2)/(1-\varrho_2)$. Then
    \begin{equation*}
    g_2=\tau_2\frac{f-\varrho_2}{f+\varrho_2}
    \end{equation*}
is twofold ramified over each of $\pm\tau_2^2,\pm1$, where $-\tau_2^2=(\sqrt{2}+1)^2>1$. By similar arguments as above, there exists an entire function $\phi_2(z)$ such that $g_2$ is written as $g_2(z)=\text{sn}(\phi_2(z))$, where $\text{sn}(\phi_2)=\text{sn}(\phi_2,-1/\tau_2^2)$. Moreover, $\phi_2'$ has order of growth strictly less than~1 when $\sigma_2(f)<1$. Thus
    \begin{equation*}
    f=-\varrho_2\frac{\text{sn}(\phi_2)+\tau_2}{\text{sn}(\phi_2)-\tau_2}.
    \end{equation*}
Iterating \eqref{yanagiharaeq13b} gives $f(z+4)^2=f(z)^2$. Therefore, we have $f(z+8)=f(z)$ and then $\text{sn}(\phi_2(z+8))=\text{sn}(\phi_2(z))$ giving $\phi_2(z+8)=\phi_2(z)+K_2$ with a period $K_2$, which is possible only when $\phi_2$ is a polynomial of degree~1 since $\phi_2'(z+8)=\phi_2'(z)$. Denote $\phi_2(z)=Cz+D$ for two constants $C,D$. Supposing that $\text{sn}(\phi_2(z_0))=1$, we then have $f(z_0)=1$. By \eqref{yanagiharaeq13b} we have $f(z_0+1)=\infty$ and it follows that $\text{sn}(\phi_2(z_0+1))=\tau_2$. On the other hand, if $\text{sn}(\phi_2(z_1))=\tau_2^2$, then $f(z_1)=-1$ and it follows from \eqref{yanagiharaeq13b} that $f(z_1+1)=\infty$ and hence $\text{sn}(\phi_2(z_1+1))=\tau_2$. Thus, $\text{sn}(\phi_2(z_0+1))=\text{sn}(\phi_2(z_1+1))$, which implies that $\phi_2(z_0)-\phi_2(z_1)=C(z_0-z_1)=\phi_2(z_0+1)-\phi_2(z_1+1)=K_2$ with a period $K_2$, but $\text{sn}(\phi_2(z_0))\not=\text{sn}(\phi_2(z_1))$, a contradiction. Therefore, \eqref{yanagiharaeq13b} cannot admit any meromorphic solutions of hyper-order strictly less than~1 when $\kappa_2$ is a constant.

From the above discussions, we obtain the following corollary which is a natural difference analogue of Steinmetz's generalization \cite{steinmetz:78} of Malmquist's~1913 result on differential equations.

\begin{corollary}\label{dtheoremc}
Let $n\in\mathbb{N}$. If the difference equation \eqref{yanagiharaeq2} with rational coefficients has a transcendental meromorphic solution of hyper-order strictly less than~1, then either $f$ satisfies \eqref{lineareq} or \eqref{driccati0} with rational coefficients or, by a transformation $f\rightarrow \alpha f$ or $f\rightarrow 1/(\alpha f)$ with an algebraic function $\alpha$ of degree at most~2, \eqref{yanagiharaeq2} reduces into one of the following equations:
   \begin{eqnarray}
   \overline{f}^2 &=& 1-f^2,\label{yanagiharaeq11 co}\\
   \overline{f}^2 &=& 1-\left(\frac{\delta f-1}{f-\delta}\right)^2,\label{yanagiharaeq12 co}\\
   \overline{f}^{2} &=& 1-\left(\frac{f+3}{f-1}\right)^{2},\label{yanagiharaeq13 co}\\
   \overline{f}^2 &=& \frac{f^2-\kappa^2}{f^2-1},\label{yanagiharaeq13a co}\\
   \overline{f}^{3} &=& 1-f^{-3},\label{yanagiharaeq15 co}
  \end{eqnarray}
where $\delta\not\equiv\pm1$ is an algebraic function of degree~2 at most and $\kappa^2\not=0,1$ is a constant.
\end{corollary}

Below we will show that equations \eqref{yanagiharaeq11 co}-\eqref{yanagiharaeq15 co} can indeed have meromorphic solutions of finite order.

Let's look at the second degree Fermat difference equations \eqref{yanagiharaeq11 co}, \eqref{yanagiharaeq12 co} and \eqref{yanagiharaeq13 co} first. For equation \eqref{yanagiharaeq11 co}, we know from \cite[Theorem~2]{Yanagihara1989difference} that the solution $f$ is represented as
$f=(\beta+\beta^{-1})/2$, where $\beta$ satisfies $\overline{\beta}=i\beta^{\pm1}$. For equation \eqref{yanagiharaeq12 co}, if we put $f=(\gamma+\gamma^{-1})/2$, then we have
    \begin{equation*}
    \frac{1}{4}\left(\overline{\gamma}-\frac{1}{\overline{\gamma}}\right)^2=-\left(\frac{\delta f-1}{f-\delta}\right)^2.
    \end{equation*}
It follows that
    \begin{equation*}
    \overline{\gamma}^2-2i\frac{\delta\gamma^2-2\gamma+\delta}{\gamma^2-2\delta\gamma+1}\overline{\gamma}-1=0.
    \end{equation*}
Solving the above equation, we get the difference Riccati equation
    \begin{equation*}
    \overline{\gamma}=\left\{-\theta\frac{(i\delta-\sqrt{1-\delta^2})\gamma+i}{\gamma-\delta+i\sqrt{1-\delta^2}}\right\}^{\theta}, \quad \theta=\pm1.
    \end{equation*}
For equation \eqref{yanagiharaeq13 co}, we have
    \begin{equation*}
    \overline{f}^{2}=-\frac{8(f+1)}{(f-1)^{2}}.
    \end{equation*}
If we put
    \begin{equation*}
    \sqrt{-8}u=\frac{\overline{f}(f-1)}{f+1}, \quad  v=\frac{1}{f+1},
    \end{equation*}
then we have
    \begin{equation*}
    u^{2}=v,  \quad  f =\frac{1}{v}-1, \quad  \overline{f}=\frac{\sqrt{-8}u}{1-2v},
    \end{equation*}
and further
    \begin{equation*}
    (\sqrt{2}\overline{u})^{2}=2\overline{v}=\frac{2}{\overline{f}+1}=\frac{2(2u^{2}-1)}{2u^{2}-\sqrt{-8}u-1}=1-\left(\frac{\sqrt{2}iu-1}{\sqrt{2}u-i}\right)^2.
    \end{equation*}
Putting $\sqrt{2}u=(\lambda+\lambda^{-1})/2$, then we have
    \begin{equation*}
    f=\frac{1-u^2}{u^2}=\frac{8\lambda^2-(\lambda^2+1)^2}{(\lambda^2+1)^2},  \quad \overline{\lambda}=\left\{-\theta\frac{-(1+\sqrt{2})\lambda+i}{\lambda-i+i\sqrt{2}}\right\}^{\theta}, \quad \theta=\pm1.
    \end{equation*}
Since the autonomous version of the difference Riccati equation is solvable explicitly in terms of exponential functions, so are equations \eqref{yanagiharaeq11 co}, \eqref{yanagiharaeq12 co} and \eqref{yanagiharaeq13 co} when $\delta$ is a constant.

Equation \eqref{yanagiharaeq13a co} can be rewritten as $\overline{f}^2f^2-(\overline{f}^2+f^2)+\kappa^2=0$, which is a special case of symmetric QRT map \cite{QuispelRobertsThompson1988,QuispelRobertsThompson1989}. By doing a suitable M\"obius transformation \cite{Baxter1982,RamaniCarsteaGrammaticosOhta2002}, for example, $f\rightarrow a(f+1)/(f-1)$, where $a$ ($|a|>1$) is a constant satisfying $2a^4-2a^2-1=0$ and $\kappa^2=a^4$, we get
     \begin{equation*}
     \overline{f}^2f^2+\overline{f}^2+f^2+4(1+4a^2)\overline{f}f+1=0,
    \end{equation*}
which is solvable in terms of Jacobi elliptic functions with finite order of growth \cite{HalburdKorhonen2007}.

For the third degree Fermat difference equation \eqref{yanagiharaeq15 co}, the solution $f$ satisfies $f^{-1}=H(\varphi)$, $\overline{f}=\eta G(\varphi)$, where $H,G$ are defined in \eqref{Weierstrassellptic0} and $\varphi$ is an entire function. Choose $\eta=1$. It follows that
    \begin{equation*}
   \frac{1+\wp'(\overline{\varphi})/\sqrt{3}}{2\wp(\overline{\varphi})}\cdot \frac{1-\wp'(\varphi)/\sqrt{3}}{2\wp(\varphi)}=1.
   \end{equation*}
Using the addition law of the Weierstrass elliptic function together with the relation $\wp'^2=4\wp^3-1$, it can be shown that this equation is solved by a polynomial of degree~1 satisfying $\overline{\varphi}=\varphi+a$, where $a$ is a constant such that $\wp'(a)=-\sqrt{3}$ and $\wp(a)=1$. It follows that the order of growth of $f$ is~2.

The following three sections contain the proof of Theorem~\ref{dtheorem}. In Section~\ref{proof1_sec}, we will first find some restrictions on the roots and degrees of the numerator and denominator of $R(z,f)$. This allows us to only consider two cases of equation \eqref{yanagiharaeq2} where $p=n,q=0$ or $p=q=n$ after a possible bilinear transformation to $f$. These two cases will be discussed in Section~\ref{proof2_sec} and Section~\ref{proof3_sec}, respectively. The results obtained in the case $p=n,q=0$ in Section~\ref{proof2_sec} are the admissible counterparts of the results on \eqref{yanagiharaeq n} by Yanagihara \cite{Yanagihara1989difference}.

\section{Restrictions on the roots and degrees of \eqref{yanagiharaeq2}}\label{proof1_sec}
We denote
    \begin{equation*}
    R(z,f) = \frac{P(z,f)}{Q(z,f)},
    \end{equation*}
where
    \begin{equation*}
    P(z,f) = a_pf^p + a_{p-1}f^{p-1}+\cdots+a_0
    \end{equation*}
and
    \begin{equation*}
    Q(z,f) = b_qf^q + b_{q-1}f^{q-1}+\cdots+b_0
    \end{equation*}
are polynomials in $f$ having no common factors and $p,q\in\mathbb{N}$. Then we have $\deg_f (P(z,f))=p$, $\deg_f (Q(z,f))=q$ and by assumption that $\deg_f (R(z,f))=\max\{p,q\}=n$. In what follows we also write
    \begin{equation}\label{P}
    P(z,f) = a_p(f-\alpha_1)^{k_1}\cdots (f-\alpha_\mu)^{k_\mu}
    \end{equation}
and
    \begin{equation}\label{Q}
    Q(z,f) = b_q(f-\beta_1)^{l_1}\cdots (f-\beta_\nu)^{l_\nu},
    \end{equation}
where the coefficients $\alpha_1,\ldots,\alpha_\mu$ and $\beta_1,\ldots,\beta_\nu$ are in general algebroid functions and $k_i$ and $l_j$ denote the order of the roots $\alpha_i$ and $\beta_j$, respectively. Without losing generality, we may suppose that the \emph{greatest common divisor} of $k_1,\ldots,k_\mu,l_1,\ldots,l_\mu$, which is denoted by $k=(k_1,\ldots,k_\mu,l_1,\ldots,l_\mu)$, is~$1$. Otherwise, after taking the $k$-th root on both sides of \eqref{yanagiharaeq2}, equation \eqref{yanagiharaeq2} reduces into a difference equation of degree~$n/k$ with meromorphic coefficients. Note that under this assumption either $P(z,f)$ or $Q(z,f)$ has at least two distinct roots when $p=0$ or $q=0$ or $p=q=n$.

Let now $f$ be an admissible meromorphic solution of \eqref{yanagiharaeq2}. For the simple case $n=1$, \eqref{yanagiharaeq2} is the linear difference equation \eqref{lineareq} or the difference Riccati equation \eqref{driccati0}. From now on, we take $n$ to be $\geq2$.
By making use of the factorizations \eqref{P} and \eqref{Q}, it follows that the roots of $P(z,f)$ and $Q(z,f)$ are $\alpha_1,\ldots,\alpha_\mu$ and $\beta_1,\ldots,\beta_\nu$, respectively. Suppose that $z_0\in\mathbb{C}$ is such that
    \begin{equation*}
    f(z_0)-\alpha_i(z_0)=0
    \end{equation*}
with multiplicity $m\in\mathbb{Z}^+$. Now, $n|mk_i$ since otherwise $z_0+1$ would be an algebraic branch point of $f$. Hence, $n\leq mk_i$ and so $m \geq n/k_i$. The same inequality holds for the roots $\beta_j$, as well. Therefore,
    \begin{equation}\label{ki}
    \overline{N}\left(r,\frac{1}{f-\alpha_i}\right) \leq \frac{k_i}{n} N\left(r,\frac{1}{f-\alpha_i}\right)
    \end{equation}
for all $i=1,\ldots,\mu$ and
    \begin{equation}\label{lj}
    \overline{N}\left(r,\frac{1}{f-\beta_j}\right) \leq \frac{l_j}{n} N\left(r,\frac{1}{f-\beta_j}\right)
    \end{equation}
for all $i=1,\ldots,\nu$. Note that the above two inequalities \eqref{ki} and \eqref{lj} still hold true for $\overline{f}-\overline{\alpha}_i$ and $\overline{f}-\overline{\beta}_j$, respectively. Now,
    \begin{equation*}
    k_1+\cdots+k_\mu \leq n, \quad l_1+\cdots+l_\nu \leq n
    \end{equation*}
and so by applying the Second Main Theorem \eqref{Yamanoi co} together with \eqref{ki} and \eqref{lj}, we have, for any $\varepsilon > 0$,
    \begin{equation*}
    \begin{split}
    (\mu+\nu-2-\varepsilon) T(r,f) &\leq \sum_{i=1}^\mu \overline{N}\left(r,\frac{1}{f-\alpha_i}\right) + \sum_{j=1}^\nu \overline{N}\left(r,\frac{1}{f-\beta_j}\right) \\
    &\leq \sum_{i=1}^\mu\frac{k_i}{n} N\left(r,\frac{1}{f-\alpha_i}\right) + \sum_{j=1}^\nu \frac{l_j}{n} N\left(r,\frac{1}{f-\beta_j}\right) \\
    &\leq \frac{1}{n}\left(\sum_{i=1}^\mu k_i + \sum_{j=1}^\nu l_j\right)T(r,f) \\
    &\leq 2T(r,f),
    \end{split}
    \end{equation*}
which implies that $\mu+\nu \leq 4$ and therefore the combined number of distinct roots of $P(z,f)$ and $Q(z,f)$ is at most~4. In particular, if $k_1+\cdots+k_\mu < n$ or $l_1+\cdots+l_\nu<n$, then $\mu+\nu\leq 3$.

Consider first the case $p=n$ and $1\leq q \leq n-1$. We show that this case cannot occur. Suppose that $N(r,f)\not=S(r,f)$. Then there are more than $S(r,f)$ points $z_0\in\mathbb{C}$ such that
    \begin{equation*}
    f(z)= C(z-z_0)^{-m} + O\left((z-z_0)^{-m+1}\right),\quad C\not=0, \quad m\in\mathbb{Z}^+
    \end{equation*}
in a neighborhood of $z_0$. Let these poles be our starting points for iteration. Note that $T(r,f(z+j))=T(r,f)+S(r,f)$, $j\in \mathbb{N}$, and all the coefficients of \eqref{yanagiharaeq2} are in $\mathcal{S}'(f)$. It follows by \cite[Lemma~3.1]{halburdrk:LMS2006} that, for an arbitrarily small $\varepsilon\geq 0$,
there are more than $S(r,f)$ points $z_0\in\mathbb{C}$ at which $f(z+1)^n$ has a pole of order at most $(1+\varepsilon)(n-q)m$ and so there are more than $S(r,f)$ poles of $f$ of order at most $(1+\varepsilon)(n-q)m/n$ at $z=z_0+1$. By continuing the above iteration it follows that there are more than $S(r,f)$ poles of $f$ of order at most $(1+\varepsilon)^2(n-q)^2m/n^2$ at $z=z_0+2$, and more than $S(r,f)$ poles of $f$ of order at most $(1+\varepsilon)^s(n-q)^{s}m/n^{s}$ at $z=z_0+s$, $s\in\mathbb{N}$.
By letting $s\to\infty$, it follows that there is necessarily a branch point of $f$ at $z_0+s_0$ for some $s_0\in\mathbb{N}$, a contradiction to our assumption that $f$ is meromorphic. Therefore, we have $N(r,f)=S(r,f)$. Also, from \eqref{yanagiharaeq2} we see that $N(r,1/(f-\beta_j))=S(r,f)$. Thus $\infty,\beta_j$, $j=1,\ldots,\nu$, are all Picard exceptional small functions of $f(z)$. By \eqref{ki} and \eqref{de relation}, we conclude that $\mu=\nu=1$. It follows that \eqref{yanagiharaeq2} assumes the following form:
    \begin{equation}\label{eqnonmero}
    \overline{f}^n=\frac{c(f-\alpha_1)^n}{(f-\beta_1)^q},
    \end{equation}
where $\alpha_1,\beta_1$ are meromorphic functions. Moreover, $\beta_1\not\equiv0$ since otherwise $\alpha_1$ is also a Picard exceptional small function of $f(z)$, which is impossible.
Denote $c_1=c^{1/n}$ and $g=(f-\beta_1)^{1/n}$. Then $c^{1/n},g^{1/n}$ are both algebroid functions of degree at most~$n$ and from the above reasoning we see that $0,\infty$ are both Picard exceptional small functions of $f(z)$. We take the $n$-th root on both sides of \eqref{eqnonmero} and then rewrite the resulting equation as follows
    \begin{equation}\label{eqnonmero ne1}
    \overline{g}^ng^q=c_1g^n-\overline{\beta}_1g^q+c_1(\beta_1-\alpha_1).
    \end{equation}
Let $u$ be an algebroid function defined by
    \begin{equation*}
    c_nu^n-\overline{\beta}_1u^{q}+c_1(\beta_1-\alpha_1)=0.
    \end{equation*}
Then $u$ is an algebroid function of degree at most $n^2$. Since $q\leq n-1$ and $\beta_1\not\equiv\alpha_1$, the above equation has at least one non-zero root $u_0$ and from \eqref{eqnonmero ne1} we see that $N(r,1/(g-u_0))=S(r,g)$, i.e., $u_0$ is also a Picard exceptional small function of $f(z)$, which is impossible. Therefore, the case where $p=n$ and $1\leq q \leq n-1$ cannot occur. This also implies that the case where~0 is a root of $P(z,f)$ of order less than~$n$ cannot occur since otherwise by doing a bilinear transformation $f\rightarrow1/f$ to \eqref{yanagiharaeq2}, we get
    \begin{equation*}
    \overline{f}^n=\frac{P_1(z,f)}{Q_1(z,f)},
    \end{equation*}
where $P_1(z,f)$ is a polynomial in $f$ of degree $n$ and $Q_1(z,f)$ is a polynomial in $f$ of degree less than $n$, which is impossible.
From the above reasoning, we conclude that $q=0$ or $q=n$ when $p=n$ and that~0 is not a root of $P(z,f)$ of order less than~$n$.

Consider now the case $q=n$ and $0\leq p\leq n-1$. Below we show that in this case~0 cannot be a root of $Q(z,f)$ of order less than~$n$. Otherwise, \eqref{yanagiharaeq2} can be written as
    \begin{equation}\label{PPrst0}
    \overline{f}^n=\frac{P(z,f)}{f^{l_0}\widehat{Q}(z,f)},
    \end{equation}
where $1\leq l_0\leq n-1$, and $\widehat{Q}(z,f)$ is a polynomial in $f$ of degree $n-l_0$. If $N(r,1/f)\not=S(r,f)$, then there are more than $S(r,f)$ points $z_0\in\mathbb{C}$ such that
    \begin{equation*}
    f(z)= C(z-z_0)^{m} + O\left((z-z_0)^{m+1}\right),\quad C\not=0,  \quad  m\in\mathbb{N}^+
    \end{equation*}
in a neighborhood of $z_0$. We now iterate \eqref{yanagiharaeq2} with these points as our starting points. Recall that $T(r,f(z+j))=T(r,f)+S(r,f)$ for all $j\in\mathbb{N}$ and all the coefficients of \eqref{yanagiharaeq2} are in $\mathcal{S}'(f)$. According to \cite[Lemma~3.1]{halburdrk:LMS2006}, for an arbitrarily small $\varepsilon\geq 0$,
by \eqref{PPrst0}, there are more than $S(r,f)$ points $z_0\in\mathbb{C}$ at which $f(z+1)^n$ has a pole of order at most $(1+\varepsilon)l_0m$
and so there are more than $S(r,f)$ poles of $f$ of order at most $(1+\varepsilon)l_0m/n$ at $z=z_0+1$. It follows that there are more than $S(r,f)$ points $z_0\in\mathbb{C}$ at which $f(z+2)^n$ has a zero of order at most $(1+\varepsilon)l_0(n-p)m$,
that is, there are more than $S(r,f)$ zeros of $f$ of order at most $(1+\varepsilon)l_0(n-p)m/n$ at $z=z_0+2$.
Then, by continuing the iteration it follows that there are more than $S(r,f)$ zeros of $f$ of order $(1+\varepsilon)^s(n-p)^sl_0^{s}m/n^{s}$ at $z=z_0+2s$, $s\in\mathbb{N}$.
By letting $s\to\infty$, it follows that there is necessarily a branch point of $f$ at $z_0+2s_0$ for some $s_0\in\mathbb{N}$, a contradiction to our assumption that $f$ is meromorphic. Therefore, $N(r,1/f)=S(r,f)$. Also, from \eqref{PPrst0} we see that $N(r,f)=S(r,f)$ since $p<q$ and then, since $l_0\leq n-1$, it follows that there exists another nonzero $\beta_j$ such that $N(r,1/(f-\beta_j)=S(r,f)$, that is to say, $f$ has at least~3 Picard exceptional small functions, which is impossible. Therefore,~$0$ cannot be a root of $Q(z,f)$ of order $<n$ when $0\leq p\leq n-1$. Now, since~0 is not a root of $P(z,f)$ either, by doing a bilinear transformation $f\rightarrow1/f$ to \eqref{yanagiharaeq2}, we get
    \begin{equation*}
    \overline{f}^n=\frac{P_2(z,f)}{Q_2(z,f)},
    \end{equation*}
where $P_2(z,f)$ and $Q_2(z,f)$ are two polynomials in $f$, both of them having degree equal to $n$. In particular, $Q_2(z,f)=f^n$ when $p=0$.

From the above discussions, we conclude that we only need to consider \eqref{yanagiharaeq2} for two cases where $p=n$, $q=0$ or $p=q=n$. Moreover, if $P(z,f)$ has two or more distinct roots, then none of them vanishes identically.

\section{Equation \eqref{yanagiharaeq2} with $p=n$ and $q=0$}\label{proof2_sec}

In Section~\ref{proof1_sec}, we have shown that the combined number of distinct roots of $P(z,f)$ is at most~3 when $p=n$ and $q=0$ and assumed that $P(z,f)$ has at least~2 distinct roots. Hence we have the following two possibilities:
    \begin{eqnarray}
    \overline{f}^n &=& c(f-\alpha_1)^{\kappa_1}(f-\alpha_2)^{\kappa_2}, \quad \alpha_1\alpha_2\not\equiv0,\label{Preq2}\\
    \overline{f}^n &=& c(f-\alpha_1)^{\mu_1}(f-\alpha_2)^{\mu_2}(f-\alpha_3)^{\mu_3}, \quad \alpha_1\alpha_2\alpha_3\not\equiv0,\label{Preq3}
    \end{eqnarray}
where $\kappa_i$, $\mu_j$ are positive integers satisfying $\kappa_1+\kappa_2=n$, $\mu_1+\mu_2+\mu_3=n$.

We first look at equation \eqref{Preq2}. Suppose that $\kappa_1+\kappa_2=n\geq3$.
Since $(\kappa_1,\kappa_2)=(n,\kappa_2)=(n,\kappa_2)=1$, we have from \eqref{ki} that
\begin{equation*}
\overline{N}\left(r,\frac{1}{f-\alpha_i}\right) \leq \frac{1}{n} N\left(r,\frac{1}{f-\alpha_i}\right), \quad i=1,2.
\end{equation*}
If $n=3$ or $n=4$, then obviously $\min\{\kappa_1,\kappa_2\}=1$; if $n\geq 5$ and $\min\{\kappa_1,\kappa_2\}\geq2$, then from \eqref{Preq2} we see that $f$ is also completely ramified at zeros, which and \eqref{de relation} imply that $2-2/n+1/2\leq 2$, a contradiction. Hence we must have $\min\{\kappa_1,\kappa_2\}=1$ and it follows that \eqref{Preq2} is of the form:
    \begin{equation}\label{Preq20before}
    \overline{f}^{n} = c(f-\alpha_1)^{n-1}(f-\alpha_2).
    \end{equation}
Moreover, since $n\geq 3$, we have $N(r,1/(f-\alpha_2))\not=S(r,f)$ for otherwise $f$ is $n$-fold ramified at $\alpha_1$-points and $(n-1)$-fold ramified at zeros, which yields a contradiction to \eqref{de relation}. Denote $c_1=[c(\alpha_2-\alpha_1)]^{1/n}$ and $c_2=1/(\alpha_2-\alpha_1)$. Then $c_1$ and $c_2$ are two algebroid functions of degree at most $n^2$ and $n$, respectively. Put
    \begin{equation}\label{Preq20}
    c_1u=\frac{\overline{f}}{f-\alpha_1}, \quad v=\frac{1}{f-\alpha_1},
    \end{equation}
where $u$ and $v$ are two algebroid functions with at most $S(r,f)$ many branch points. By the First Main Theorem, we have $T(r,v)=T(r,f)+S(r,f)$. Thus all the coefficients of \eqref{Preq20before} are small functions of $v$. Then we have
    \begin{equation*}
    f=\frac{1}{v}+\alpha_1, \quad \overline{f}=\frac{c_1u}{v}
    \end{equation*}
and from \eqref{Preq20before} and \eqref{Preq20} that
    \begin{equation}\label{Preq22}
    u^{n}=c_2-v.
    \end{equation}
It also follows that
    \begin{equation}\label{Preq23}
    \frac{c_1u-\overline{\alpha}_2v}{v}=\overline{f}-\overline{\alpha}_2.
    \end{equation}
Since zeros of $\overline{f}-\overline{\alpha}_2$ are $n$-fold ramified, we see that zeros of $F=c_1u-\overline{\alpha}_2v$ are also $n$-fold ramified. Let $u_0$ and $v_0$ be a pair of roots of the following equations
    \begin{equation*}
    u_0^{n}=c_2-v_0, \quad c_1u_0-\overline{\alpha}_2v_0=0.
    \end{equation*}
Eliminating $u_0$ from the above two equations gives
    \begin{equation}\label{PPreq6}
    \overline{\alpha}_2^{n}v_0^{n}+c_1^{n}v_0-c_1^{n}c_2=0,
    \end{equation}
from which we see that $v_0$ is an algebroid function of degree at most $n^2$. Moreover, it is easy to see that $c_2$ is not a root of \eqref{PPreq6}. If \eqref{PPreq6} has only one root, then, since $n\geq3$, by Vieta's formulas we see that $v_0=0$ and it follows from \eqref{PPreq6} that $c_1^{n}c_2\equiv0$, a contradiction.
Therefore, \eqref{PPreq6} has at least two distinct roots, say, $v_1$ and $v_2$. Rewrite \eqref{Preq23} as
    \begin{equation}\label{Preq23hh}
    -\frac{\overline{\alpha}_2^{n}v^n+c_1^nv-c_1^nc_2}{v^n}=\overline{f}^n-\overline{\alpha}_2^{n}.
    \end{equation}
Denote the order of the roots $v_1$ and $v_2$ by $l_1$ and $l_2$, respectively. Then $1\leq l_1,l_2\leq n-1$ and we have
     \begin{equation*}
    \overline{\alpha}_2^{n}v^n+c_1^nv-c_1^nc_2=\overline{\alpha}_2^{n}(v-v_1)^{l_1}(v-v_2)^{l_2}P(z,v),
    \end{equation*}
where $P(z,v)$ is a polynomial in $v$ of degree $n-l_1-l_2$ with small functions of $v$ as coefficients.
Suppose that $z_0\in\mathbb{C}$ is such that
    \begin{equation*}
    f(z_0+1)-\alpha_2(z_0+1)=0
    \end{equation*}
with multiplicity $n_0\geq{n}$. If $v(z_0)-v_1(z_0)=0$ with multiplicity $m\in{\mathbb{Z}^+}$, then from \eqref{Preq23hh}, we have $n_0|ml_1$. Hence, $n_0\leq ml_1$ and so $m \geq 2$. Therefore,
    \begin{equation*}
    \overline{N}\left(r,\frac{1}{v-v_1}\right) \leq \frac{1}{2} N\left(r,\frac{1}{v-v_1}\right).
    \end{equation*}
The same inequality holds for $v_2$, as well. Thus $v$ is completely ramified at $v_1$-points and $v_2$-points. But from \eqref{Preq22} we see that $v$ is $n$-fold ramified at poles and $c_2$-points, which yields a contradiction to \eqref{de relation} since $n\geq 3$. Hence $n=2$ and we therefore have
    \begin{equation}\label{PPreq8}
    \overline{f}^{2} = c(f-\alpha_1)(f-\alpha_2).
    \end{equation}

In what follows, we will still use the expressions in \eqref{Preq20} to transform \eqref{PPreq8} into \eqref{Preq22} for the case $n=2$. Moreover, from the above reasoning we know that
    \begin{equation}\label{PPreq9}
    \overline{\alpha}_2^{2}v_0^{2}+c_1^{2}v_0-c_1^{2}c_2=0,
    \end{equation}
which has two distinct roots if
    \begin{equation}\label{PPreq9ttt}
    (4\overline{\alpha}_2^2c_2+c_1^2)c_1^2\not=0.
    \end{equation}
If \eqref{PPreq9ttt} holds, then $v$ must be completely ramified at poles, $c_2$-points, $v_1$-points and $v_2$-points; if not, then we have
    \begin{equation}\label{PPreq9tttg0}
    c_1^2+4\overline{\alpha}_2^2c_2=\frac{4\overline{\alpha}_2^2+c(\alpha_2-\alpha_1)^2}{\alpha_2-\alpha_1}=0.
    \end{equation}
Rewrite \eqref{PPreq8} as
    \begin{equation*}
    \overline{f}^{2} = c(g^2-c_4^2),
    \end{equation*}
where $g=f-c_3$ is an algebroid function with at most $S(r,f)$ many branch points and $c_3=(\alpha_1+\alpha_2)/2$, $c_4=(\alpha_2-\alpha_1)/2=1/(2c_2)$ are two algebroid functions of degree at most~2. Put
    \begin{equation*}
    g=\frac{c_4}{2}\left(\beta+\frac{1}{\beta}\right).
    \end{equation*}
Then substitution yields
    \begin{equation*}
    \overline{f}^{2} =\frac{cc_4^2}{4}\left(\beta-\frac{1}{\beta}\right)^2,
    \end{equation*}
and, without loss of generality, we may suppose that
    \begin{equation}\label{PPreq10}
    \beta^2-\frac{2}{c^{1/2}c_4}\left(\overline{\alpha}_1+\frac{1}{\overline{v}}\right)\beta-1=0.
    \end{equation}
Note that $4c_2^2c_4^2=1$ and $c=c_1^2c_2$. We have $4cc_4^2=c_1^2/c_2$. The discriminant of the quadratic equation \eqref{PPreq10} with respect to $\beta$,
    \begin{equation}\label{PPreq11before}
    \frac{4(\overline{\alpha}_1\overline{v}+1)^2+4cc_4^2\overline{v}^2}{cc_4^2\overline{v}^2}=\frac{4c_2(\overline{\alpha}_1\overline{v}+1)^2+c_1^2\overline{v}^2}{c_2cc_4^2\overline{v}^2}=\frac{4\overline{f}^2}{cc_4^2}+4=\frac{4(cc_4^2+\overline{f}^2)}{cc_4^2}=\frac{4g^2}{c_4^2},
    \end{equation}
implies that $\beta$ is meromorphic apart from at most $S(r,f)$ many points.
Moreover, roots of the equation
    \begin{equation}\label{PPreq11}
    (4\overline{\alpha}_1^2c_2+c_1^2)\overline{v}^2+8c_2\overline{\alpha}_1\overline{v}+4c_2=0
    \end{equation}
with respect to $\overline{v}$ must be either $c_2$ or coincide with $v_1,v_2$.

If $c_2$ is a root of \eqref{PPreq11}, then we have
    \begin{equation*}
     \frac{4\overline{\alpha}_1^2c_2+c_1^2}{(\overline{\alpha}_2-\overline{\alpha}_1)^2}+\frac{8c_2\overline{\alpha}_1}{\overline{\alpha}_2-\overline{\alpha}_1}+4c_2=0,
    \end{equation*}
which yields \eqref{PPreq9tttg0}. On the other hand, by exchanging the roles of $\alpha_1$ and $\alpha_2$, we also have
    \begin{equation}\label{PPreq9tttg0 ex}
    c_1^2+4\overline{\alpha}_1^2c_2=0.
    \end{equation}
From \eqref{PPreq9tttg0} and \eqref{PPreq9tttg0 ex} we see that $\alpha_2=-\alpha_1$ and it follows that $c=-\overline{\alpha}_1^2/\alpha_1^2$. By doing the transformation $f\rightarrow\alpha_1f$, then \eqref{PPreq8} becomes
    \begin{equation*}
    \overline{f}^2=1-f^2,
    \end{equation*}
which is the equation \eqref{yanagiharaeq11}.

If $c_2$ is not a root of \eqref{PPreq11}, then roots of \eqref{PPreq11} must coincide with $v_1,v_2$,
which means that
    \begin{equation}\label{PPreq12aabefore e}
    \frac{\overline{\overline{\alpha}}_2^{2}}{4\overline{\alpha}_1^2c_2+c_1^2}=\frac{\overline{c}_1^{2}}{8c_2\overline{\alpha}_1}=\frac{-\overline{c}_1^{2}\overline{c}_2}{4c_2}.
   \end{equation}
Recall that $c_1=[c(\alpha_2-\alpha_1)]^{1/2}$ and $c_2=1/(\alpha_2-\alpha_1)$. The second equation of \eqref{PPreq12aabefore e} gives $-2\overline{c}_2\overline{\alpha}_1=1$, i.e., $\overline{\alpha}_2=-\overline{\alpha}_1$. It follows from the first equation of \eqref{PPreq12aabefore e} that
    \begin{equation*}\label{PPreq12aabefore}
    \alpha_1^2\overline{c}c+\overline{\alpha}_1^2\overline{c}+\overline{\overline{\alpha}}_1^2=0.
     \end{equation*}
Denote $d=c\alpha_1^2/\overline{\alpha}_1^2$. Then $d$ is a meromorphic function
satisfying $\overline{d}(d+1)+1=0$. By doing the transformation $f\rightarrow\alpha_1f$, then \eqref{PPreq8} becomes
    \begin{equation*}\label{PPreq12aa}
    \overline{f}^2=d(f^2-1),
    \end{equation*}
which is the equation \eqref{yanagiharaeq11a}.

We now look at equation \eqref{Preq3}. Denote $k_1=(n,\mu_1)$, $k_2=(n,\mu_2)$, $k_3=(n,\mu_3)$ and let $K_1=n/k_1$, $K_2=n/k_2$, $K_3=n/k_3$. Obviously, $K_i\geq2$, $i=1,2,3$, and $1/K_1+1/K_2+1/K_3=(k_1+k_2+k_3)/n\leq 1$. Now $f$ is $K_i$-fold ramified at $\alpha_i$-points, $i=1,2,3$, which and \eqref{de relation} imply that
    \begin{equation*}
1-\frac{1}{K_1}+1-\frac{1}{K_2}+1-\frac{1}{K_3}\leq 2,
    \end{equation*}
and so $1/K_1+1/K_2+1/K_3=(k_1+k_2+k_3)/n=1$. Thus we have the following three possibilities:
    \begin{equation*}
    \begin{split}
    1,  \quad  & K_1=2, K_2=3, K_3=6;\\
    2,  \quad  & K_1=2, K_2=4, K_3=4;\\
    3,  \quad  & K_1=3, K_2=3, K_3=3,
    \end{split}
    \end{equation*}
which implies that \eqref{Preq3} assumes one of the following three forms:
    \begin{eqnarray}
    \overline{f}^3 &=& c(f-\alpha_1)(f-\alpha_2)(f-\alpha_3),  \quad  \alpha_1\alpha_2\alpha_3\not\equiv0,\label{PQreq24}\\
    \overline{f}^4 &=& c(f-\alpha_1)^{2}(f-\alpha_2)(f-\alpha_3),  \quad \alpha_1\alpha_2\alpha_3\not\equiv0,\label{PQreq25}\\
    \overline{f}^6 &=& c(f-\alpha_1)^{3}(f-\alpha_2)^{2}(f-\alpha_3),  \quad \alpha_1\alpha_2\alpha_3\not\equiv0.\label{PQreq26}
   \end{eqnarray}

Let's first consider equation \eqref{PQreq24}. Put
    \begin{equation}\label{PQreq27before}
    \frac{\overline{f}}{f-\alpha_1}=c_1u, \quad \frac{1}{f-\alpha_1}=b_2v+b_1,
   \end{equation}
where $u$ and $v$ are two algebroid functions with at most $S(r,f)$ many branch points and
    \begin{equation}\label{coefficients1}
    \begin{split}
     c_1 &= [4cb_2^2(\alpha_2-\alpha_1)(\alpha_3-\alpha_1)]^{1/3},\\
     b_1 &= \frac{1}{2}\left(\frac{1}{\alpha_2-\alpha_1}+\frac{1}{\alpha_3-\alpha_1}\right),\\
     b_2 &= \frac{1}{2}\left(\frac{1}{\alpha_2-\alpha_1}-\frac{1}{\alpha_3-\alpha_1}\right)
     \end{split}
   \end{equation}
are algebroid functions of degrees at most~$9,3,3$, respectively.
By the First Main Theorem, we have $T(r,v)=T(r,f)+S(r,f)$. Thus all the coefficients of \eqref{PQreq24} are small functions of $v$. Then we have
    \begin{equation*}
    f=\alpha_1+\frac{1}{b_2v+b_1}, \quad \overline{f}=\frac{c_1u}{b_2v+b_1}
   \end{equation*}
and from \eqref{PQreq24} and \eqref{PQreq27before} that
    \begin{equation}\label{PQreq27}
    v^2=4u^3+1.
   \end{equation}
It also follows that
    \begin{equation}\label{PQreq28ssdw}
   \frac{c_1u-\overline{\alpha}_2(b_2v+b_1)}{b_2v+b_1}=\overline{f}-\overline{\alpha}_2.
   \end{equation}
Since zeros of $\overline{f}-\overline{\alpha}_2$ are threefold ramified, we see that zeros of $F=c_1u-\overline{\alpha}_2(b_2v+b_1)$
are also threefold ramified. Let $u_0$ and $v_0$ be a pair of roots of the following equations
    \begin{equation*}
    v_0^2=4u_0^3+1, \quad c_1u_0-\overline{\alpha}_2(b_2v_0+b_1)=0.
    \end{equation*}
Eliminating $u_0$ from the above two equations gives
    \begin{equation}\label{PPreq29}
    4\overline{\alpha}_2^3(b_2v_0+b_1)^3-c_1^3(v_0^2-1)=0,
    \end{equation}
from which we see that $v_0$ is an algebroid function of degree at most~$9$.
Moreover, it is seen that~$\pm1$ are not roots of \eqref{PPreq29} for otherwise we have $b_2{\pm}b_1=0$, a contradiction. Rewrite \eqref{PQreq28ssdw} as
   \begin{equation}\label{PQreq28after}
   -\frac{4\overline{\alpha}_2^3(b_2v+b_1)^3-c_1^3(v^2-1)}{4(b_2v+b_1)^3}=\overline{f}^3-\overline{\alpha}_2^3.
   \end{equation}
From equation \eqref{PQreq27}, we see that $v$ is threefold ramified over~$\pm1,\infty$.
This implies that \eqref{PPreq29} can admit only one root for otherwise by a similar analysis as to \eqref{Preq23hh} we obtain from \eqref{PQreq28after} that \eqref{PPreq29} has at least two roots, say, $v_1$ and $v_2$, distinct from~$\pm1,\infty$ such that $v$ is completely ramified at $v_1$-points and $v_2$-points, which is impossible by Theorem~\ref{completelyrm}. Expanding \eqref{PPreq29} gives
    \begin{equation*}
    4\overline{\alpha}_2^3b_2^3v_0^3+(12\overline{\alpha}_2^3b_2^2b_1-c_1^3)v_0^2+12\overline{\alpha}_2^3b_2b_1^2v_0+4\overline{\alpha}_2^3b_1^3+c_1^3=0,
    \end{equation*}
and so by Vieta's formulas, we get
    \begin{equation}\label{PPreq29extras}
    \frac{12\overline{\alpha}_2^3b_2^2b_1-c_1^3}{4\overline{\alpha}_2^3b_2^3}=-3v_0, \quad \frac{3b_1^2}{b_2^2}=3v_0^2.
    \end{equation}
From the second equation of \eqref{PPreq29extras} we get $v_0=\pm b_1/b_2$. If $v_0=-b_1/b_2$, then substitution into \eqref{PPreq29} yields $b_1^2=b_2^2$, a contradiction. Hence $v_0=b_1/b_2$. Combining this equation with the first equation of \eqref{PPreq29extras}, we get
    \begin{equation*}
    c_1^3=24\overline{\alpha}_2^3b_2^2b_1.
    \end{equation*}
Similarly, the above equation also holds when replacing $\overline{\alpha}_2$ in \eqref{PQreq28ssdw} with $\overline{\alpha}_1$ and $\overline{\alpha}_3$,
which implies that $\alpha_1^3=\alpha_2^3=\alpha_3^3$. Thus, $\alpha_2=\eta\alpha_1$, $\alpha_3=\eta^2\alpha_1$, where $\eta$ is the cubic root of~1 satisfying $\eta^2+\eta+1=0$.
Moreover, we have
     \begin{equation*}
    c_1^3=4cb_2^2(\alpha_2-\alpha_1)(\alpha_3-\alpha_1)=24\overline{\alpha}_1^3b_2^2b_1,
     \end{equation*}
which yields $c=-\overline{\alpha}_1^3/\alpha_1^3$.
By doing the transformation $f\rightarrow\alpha_1f$, then \eqref{PQreq24} becomes
    \begin{equation*}\label{PQreq38}
     \overline{f}^3=1-f^3,
    \end{equation*}
which is the equation \eqref{yanagiharaeq14}.

We now consider equation \eqref{PQreq25}. Put
    \begin{equation}\label{PQreq39before}
    \frac{\overline{f}}{f-\alpha_1}=c_1u, \quad \frac{1}{f-\alpha_1}=b_2v+b_1,
    \end{equation}
where $u$ and $v$ are two algebroid functions with at most $S(r,f)$ many branch points and
    \begin{equation}\label{coefficients2}
    \begin{split}
    c_1&=[-b_2^2c(\alpha_2-\alpha_1)(\alpha_3-\alpha_1)]^{1/4},\\
    b_1&=\frac{1}{2}\left(\frac{1}{\alpha_2-\alpha_1}+\frac{1}{\alpha_3-\alpha_1}\right),\\
    b_2&=\frac{1}{2}\left(\frac{1}{\alpha_2-\alpha_1}-\frac{1}{\alpha_3-\alpha_1}\right)
    \end{split}
   \end{equation}
are algebroid functions of degrees at most~$16,4,4$, respectively.
By the First Main Theorem, we have $T(r,v)=T(r,f)+S(r,f)$. Thus all the coefficients of \eqref{PQreq25} are small functions of $v$. Then we have
    \begin{equation*}
    f=\frac{1}{b_2v+b_1}+\alpha_1, \quad \overline{f}=\frac{c_1u}{b_2v+b_1}
   \end{equation*}
and from \eqref{PQreq25} and \eqref{PQreq39before} that
    \begin{equation}\label{PQreq39}
    v^2=1-u^4.
   \end{equation}
It also follows that
   \begin{equation}\label{PQreq39ted}
   \frac{c_1u-\overline{\alpha}_2(b_2v+b_1)}{b_2v+b_1}=\overline{f}-\overline{\alpha}_2.
   \end{equation}
Since zeros of $\overline{f}-\overline{\alpha}_2$ are fourfold ramified, we see that zeros of $F=c_1u-\overline{\alpha}_2(b_2v+b_1)$ are also fourfold ramified. Let $u_0$ and $v_0$ be a pair of roots of the following equations
    \begin{equation*}
    v_0^2=1-u_0^4, \quad c_1u_0-\overline{\alpha}_2b_2v_0-\overline{\alpha}_2b_1=0.
    \end{equation*}
Eliminating $u_0$ from the above two equations gives
    \begin{equation}\label{PQreq41}
    \overline{\alpha}_2^4(b_2v_0+b_1)^4-c_1^4(1-v_0^2)=0,
    \end{equation}
from which we see that $v_0$ is an algebroid function of degree at most~$16$.
Moreover, it is seen that~$\pm1$ are not roots of \eqref{PQreq41} for otherwise we have $b_2{\pm}b_1=0$, a contradiction. Rewrite \eqref{PQreq39ted} as
   \begin{equation}\label{PQreq39ted'}
   \frac{c_1^4(1-v^2)-\overline{\alpha}_2^4(b_2v+b_1)^4}{(b_2v+b_1)^4}=\overline{f}^4-\overline{\alpha}_2^4.
   \end{equation}
From equation \eqref{PQreq39}, we see that $v$ is completely ramified over~$\pm1,\infty$. This implies that \eqref{PQreq41} can admit only one root for otherwise by a similar analysis as to \eqref{Preq23hh} we obtain from \eqref{PQreq39ted'} that \eqref{PQreq41} has at least two roots, say, $v_1$ and $v_2$, distinct from~$\pm1,\infty$ such that $v$ is completely ramified at $v_1$-points and $v_2$-points, which is impossible by Theorem~\ref{completelyrm}. Expanding \eqref{PQreq41} gives
    \begin{equation*}
    \overline{\alpha}_2^4(b_2^4v_0^4+4b_2^3b_1v_0^3+\ldots+b_1^4)-c_1^4(1-v_0^2)=0,
    \end{equation*}
and so by Vieta's formulas, we obtain $4b_1/b_2=-4v_0$, which gives $b_2v_0+b_1=0$ and it follows from \eqref{PQreq41} that $v_0^2=1$, a contradiction. Therefore, \eqref{PQreq25} cannot admit any meromorphic solutions.

Finally, let's consider equation \eqref{PQreq26}. Put
    \begin{equation}\label{PQreq51before}
    \frac{\overline{f}}{f-\alpha_1}=c_1u, \quad \frac{1}{f-\alpha_1}=b_2v-b_1,
    \end{equation}
where $u$ and $v$ are two algebroid functions with at most $S(r,f)$ many branch points and
    \begin{equation}\label{coefficients3}
    \begin{split}
    c_1&=[-c^{6/n}b_2^3(\alpha_2-\alpha_1)^2(\alpha_3-\alpha_1)]^{1/6},\\
    b_1&=-\frac{1}{3}\left(\frac{2}{\alpha_2-\alpha_1}+\frac{1}{\alpha_3-\alpha_1}\right),\\
    b_2&=\frac{1}{\alpha_2-\alpha_1}+b_1=\frac{1}{3}\left(\frac{1}{\alpha_2-\alpha_1}-\frac{1}{\alpha_3-\alpha_1}\right)
    \end{split}
    \end{equation}
are algebroid functions of degrees at most~$36,6,6$, respectively.
By the First Main Theorem, we have $T(r,v)=T(r,f)+S(r,f)$. Thus all the coefficients of \eqref{PQreq26} are small functions of $v$. Then we have
    \begin{equation*}
    f=\frac{1}{b_2v-b_1}+\alpha_1, \quad \overline{f}=\frac{c_1u}{b_2v-b_1}
   \end{equation*}
and from \eqref{PQreq26} and \eqref{PQreq51before} that
    \begin{equation}\label{PQreq51}
    u^6=(v-1)^2(v+2)=v^3-3v+2.
   \end{equation}
It also follows that
   \begin{equation}\label{PQreq51des}
   \frac{c_1u-\overline{\alpha}_3(b_2v-b_1)}{b_2v-b_1}=\overline{f}-\overline{\alpha}_3.
   \end{equation}
Since zeros of $\overline{f}-\overline{\alpha}_3$ are sixfold ramified, we see that zeros of $F=c_1u-\overline{\alpha}_3(b_2v-b_1)$ are also sixfold ramified. Let $u_0$ and $v_0$ be a pair of roots of the following equations
    \begin{equation*}
    u_0^6=v_0^3-3v_0+2, \quad c_1u_0-\overline{\alpha}_3b_2v_0+\overline{\alpha}_3b_1=0.
    \end{equation*}
Eliminating $u_0$ from the above two equations gives
    \begin{equation}\label{PQreq53}
    \overline{\alpha}_2^6(b_2v_0-b_1)^6-c_1^6(v_0^3-3v_0+2)=0,
    \end{equation}
from which we see that $v_0$ is an algebroid function of degree at most~$36$.
Moreover, it is seen that~$1$ and~$-2$ are not roots of \eqref{PQreq53} for otherwise we have $b_2-b_1=0$ or $2b_2+b_1=0$, a contradiction. Rewrite \eqref{PQreq51des} as
   \begin{equation}\label{PQreq51des'}
   \frac{c_1^6(v^3-3v+2)-\overline{\alpha}_3^6(b_2v-b_1)^6}{(b_2v-b_1)^6}=\overline{f}^6-\overline{\alpha}_3^6.
   \end{equation}
From equation \eqref{PQreq51}, we see that $v$ is completely ramified over~$1,-2,\infty$. This implies that \eqref{PQreq53} can admit only one root for otherwise by a similar analysis as to \eqref{Preq23hh} we obtain from \eqref{PQreq51des'} that \eqref{PQreq53} has at least two roots, say, $v_1$ and $v_2$, distinct from~$1,-2,\infty$ such that $v$ is completely ramified at $v_1$-points and $v_2$-points, which is impossible by Theorem~\ref{completelyrm}. Expanding \eqref{PQreq53} gives
    \begin{equation*}
    \overline{\alpha}_2^6(b_2^6v_0^6-6b_2^5b_1v_0^5+\ldots+b_1^6)-c_1^6(v_0^3-3v_0+2)=0,
    \end{equation*}
and so by Vieta's formulas, we obtain $-6b_1/b_2=-6v_0$, which gives $b_2v_0-b_1=0$ and it follows from \eqref{PQreq53} that $v_0^3-3v_0+2=0$, a contradiction. Therefore, \eqref{PQreq26} cannot admit any meromorphic solutions.

\section{Equation \eqref{yanagiharaeq2} with $p=q=n$}\label{proof3_sec}

In this section, we consider the cases where $Q(z,f)$ has one, two or three roots.

\subsection{$Q(z,f)$ has only one root}\label{proof31_sec}
In Section~\ref{proof1_sec}, we have shown that the combined number of distinct roots of $P(z,f)$ and $Q(z,f)$ is at most~$4$ and assumed that $P(z,f)$ has at least~2 distinct roots. Hence we have the following two possibilities:
     \begin{eqnarray}
    \overline{f}^n &=& \frac{c(f-\alpha_1)^{\kappa_1}(f-\alpha_2)^{\kappa_2}}{(f-\beta_1)^n}, \quad \alpha_1\alpha_2\not\equiv0,\label{Preq2s}\\
    \overline{f}^n &=& \frac{c(f-\alpha_1)^{\mu_1}(f-\alpha_2)^{\mu_2}(f-\alpha_3)^{\mu_3}}{(f-\beta_1)^n}, \quad \alpha_1\alpha_2\alpha_3\not\equiv0,\label{Preq3s}
     \end{eqnarray}
where $\kappa_i$, $\mu_j$ are positive integers satisfying $\kappa_1+\kappa_2=n$, $\mu_1+\mu_2+\mu_3=n$.

We now look at equation \eqref{Preq2s}. Suppose $\kappa_1+\kappa_2=n\geq3$. By a similar analysis on the $\alpha_i$-points of $f$ in \eqref{Preq2s} as to the $\alpha_i$-points of $f$ in \eqref{Preq2} in Section~\ref{proof2_sec}, we can obtain that \eqref{Preq2s} is of the form:
    \begin{equation}\label{Preq20sreqa}
    \overline{f}^{n} = \frac{c(f-\alpha_1)^{n-1}(f-\alpha_2)}{(f-\beta_1)^n}.
    \end{equation}
Moreover, we have $N(r,1/(f-\alpha_2))\not=S(r,f)$. Denote $c_1=[c(\alpha_2-\alpha_1)]^{1/n}$ and $c_2=1/(\alpha_2-\alpha_1)$. Then $c_1$ and $c_2$ are two algebroid functions of degree at most $n^2$ and $n$, respectively. Put
    \begin{equation}\label{Preq20s}
    c_1u=\frac{\overline{f}(f-\beta_1)}{f-\alpha_1}, \quad v=\frac{1}{f-\alpha_1},
    \end{equation}
where $u$ and $v$ are two algebroid functions with at most $S(r,f)$ many branch points. By the First Main Theorem, we have $T(r,v)=T(r,f)+S(r,f)$. Thus all the coefficients of \eqref{Preq20sreqa} are small functions of $v$. Then we have
    \begin{equation*}\label{Preq21s}
    f =\frac{1}{v}+\alpha_1, \quad \overline{f}=\frac{c_1u}{1+(\alpha_1-\beta_1)v}
    \end{equation*}
and from \eqref{Preq20sreqa} and \eqref{Preq20s} that
    \begin{equation}\label{Preq22s}
    u^{n}=c_2-v.
    \end{equation}
It also follows that
    \begin{equation}\label{Preq23s}
    \frac{c_1u-\overline{\alpha}_2[1+(\alpha_1-\beta_1)v]}{1+(\alpha_1-\beta_1)v}=\overline{f}-\overline{\alpha}_2.
    \end{equation}
Since zeros of $\overline{f}-\overline{\alpha}_2$ are $n$-fold ramified, we see that zeros of $F=c_1u-\overline{\alpha}_2[(\alpha_1-\beta_1)v+1]$ are also $n$-fold ramified. Let $u_0$ and $v_0$ be a pair of roots of the following equations
    \begin{equation*}
    u_0^{n}=c_2-v_0, \quad c_1u_0-\overline{\alpha}_2[(\alpha_1-\beta_1)v_0+1]=0.
    \end{equation*}
Eliminating $u_0$ from the above two equations gives
    \begin{equation*}
    \overline{\alpha}_2^{n}[(\alpha_1-\beta_1)v_0+1]^{n}+c_1^{n}v_0-c_1^{n}c_2=0,
    \end{equation*}
i.e.,
     \begin{equation}\label{PPreq6s}
    \overline{\alpha}_2^{n}[(\alpha_1-\beta_1)^nv_0^n+n(\alpha_1-\beta_1)^{n-1}v_0^{n-1}+\ldots+1]+c_1^{n}v_0-c_1^{n}c_2=0,
    \end{equation}
from which we see that $v_0$ is an algebroid function of degree at most $n^2$. Moreover, it is seen that $c_2$ is not a root of \eqref{PPreq6s}. If \eqref{PPreq6s} has only one root, then, since $n\geq3$, by Vieta's formulas, we see that $v_0=-1/(\alpha_1-\beta_1)$ and it follows from \eqref{PPreq6s} that $-1/(\alpha_1-\beta_1)=c_2$, which gives $\alpha_2=\beta_1$, a contradiction. Therefore, \eqref{PPreq6s} has at least two distinct roots, say, $v_1$ and $v_2$. Rewrite \eqref{Preq23s} as
    \begin{equation}\label{Preq23s'}
    -\frac{\overline{\alpha}_2^n[1+(\alpha_1-\beta_1)v]^n+c_1^{n}v-c_1^{n}c_2}{[1+(\alpha_1-\beta_1)v]^n}=\overline{f}^n-\overline{\alpha}_2^n.
    \end{equation}
By a similar analysis as to \eqref{Preq23hh} in Section~\ref{proof2_sec}, we obtain from \eqref{Preq23s'} that $v$ is completely ramified at $v_1$-points and $v_2$-points. Moreover, from \eqref{Preq22s}, we see that $v$ is $n$-fold ramified at poles and $c_2$-points, which is impossible by \eqref{de relation} since $n\geq3$. Hence $n=2$ and we therefore have
    \begin{equation}\label{PPreq8s}
    \overline{f}^{2} = \frac{c(f-\alpha_1)(f-\alpha_2)}{(f-\beta_1)^2}.
    \end{equation}

In what follows, we will still use the expressions in \eqref{Preq20s} to transform \eqref{PPreq8s} into \eqref{Preq22s} for the case $n=2$. Moreover, we have
    \begin{equation}\label{PPreq9sd}
    \overline{\alpha}_2^{2}[(\alpha_1-\beta_1)^2v_0^2+2(\alpha_1-\beta_1)v_0+1]+c_1^{2}v_0-c_1^{2}c_2=0,
    \end{equation}
which has two distinct roots if
    \begin{equation*}
    \begin{split}
    [2\overline{\alpha}_2^{2}&(\alpha_1-\beta_1)+c_1^{2}]^2-4\overline{\alpha}_2^{2}(\alpha_1-\beta_1)^2(\overline{\alpha}_2^{2}-c_1^{2}c_2)\not\equiv0.
    \end{split}
    \end{equation*}
Recalling that $c_1=[c(\alpha_2-\alpha_1)]^{1/2}$ and $c_2=1/(\alpha_2-\alpha_1)$, we have
    \begin{equation}\label{PPreq9sddd'}
    \begin{split}
    [4\overline{\alpha}_2^{2}(\alpha_1-\beta_1)(\alpha_2-\beta_1)+c(\alpha_2-\alpha_1)^{2}]c_1^{2}c_2\not\equiv0.
    \end{split}
    \end{equation}
If \eqref{PPreq9sddd'} holds, then $v$ must be completely ramified at poles, $c_2$-points, $v_1$-points and $v_2$-points; if not, then we have
    \begin{equation}\label{PPreq11srt}
    4\overline{\alpha}_2^{2}(\alpha_1-\beta_1)(\alpha_2-\beta_1)+c(\alpha_2-\alpha_1)^{2}=0.
    \end{equation}
Rewrite \eqref{PPreq8s} as
    \begin{equation}\label{PPreq8sextra}
    \overline{f}^{2}(f-\beta_1)^2 = c(g^2-c_4^2),
    \end{equation}
where $g=f-c_3$ is an algebroid function with at most $S(r,f)$ many branch points and $c_3=(\alpha_1+\alpha_2)/2$, $c_4=(\alpha_2-\alpha_1)/2=1/(2c_2)$ are two algebroid functions of degree both at most~2. Put
    \begin{equation*}
    g=\frac{c_4}{2}\left(\beta+\frac{1}{\beta}\right).
    \end{equation*}
Then substitution yields
    \begin{equation*}
    \overline{f}^{2}(f-\beta_1)^2 =\frac{cc_4^2}{4}\left(\beta-\frac{1}{\beta}\right)^2,
    \end{equation*}
and, without loss of generality, we may suppose that
    \begin{equation*}
    \frac{c^{1/2}c_4}{2}\left(\beta-\frac{1}{\beta}\right)=\left(\overline{\alpha}_1+\frac{1}{\overline{v}}\right)\left[\frac{c_4}{2}\left(\beta+\frac{1}{\beta}\right)+c_3-\beta_1\right].
    \end{equation*}
i.e.,
    \begin{equation}\label{PPreq10s}
    \left[1-\frac{1}{c^{1/2}}\left(\overline{\alpha}_1+\frac{1}{\overline{v}}\right)\right]\beta^2-\frac{2(c_3-\beta_1)}{c^{1/2}c_4}\left(\overline{\alpha}_1+\frac{1}{\overline{v}}\right)\beta-\left[1+\frac{1}{c^{1/2}}\left(\overline{\alpha}_1+\frac{1}{\overline{v}}\right)\right]=0.
    \end{equation}
Recall that $f=\alpha_1+1/v$ and $f=g+c_3$. From \eqref{PPreq8sextra}, we have $\overline{f}^2=c(g^2-c_4^2)/(f-\beta_1)^2$. Then the discriminant of the quadratic equation \eqref{PPreq10s} with respect to $\beta$,
    \begin{equation}\label{PPreq10stur}
    \begin{split}
     \frac{4[(c_3-\beta_1)^2-c_4^2]}{cc_4^2}\cdot\left(\overline{\alpha}_1+\frac{1}{\overline{v}}\right)^2+4=\frac{4[(c_3-\beta_1)^2-c_4^2](g^2-c_4^2)}{c_4^2(g+c_3-\beta_1)^2}+4=\frac{4[(c_3-\beta_1)g+c_4^2]^2}{c_4^2(g+c_3-\beta_1)^2},
    \end{split}
    \end{equation}
implies that $\beta$ is meromorphic apart from at most $S(r,f)$ many points.
Moreover, roots of the equation
     \begin{equation}\label{PPreq11s}
     \begin{split}
    [4(\alpha_1-\beta_1)(\alpha_2-\beta_1)\overline{\alpha}_1^2+c(\alpha_2-\alpha_1)^2]\overline{v}^2+8(\alpha_1-\beta_1)(\alpha_2-\beta_1)\overline{\alpha}_1\overline{v}+4(\alpha_1-\beta_1)(\alpha_2-\beta_1)=0,
    \end{split}
    \end{equation}
with respect to $\overline{v}$ must be $c_2$ or coincide with $v_1,v_2$.

If $c_2=1/(\alpha_2-\alpha_1)$ is a root of \eqref{PPreq11s}, then we have
    \begin{equation*}
    \begin{split}
    4(\alpha_1-\beta_1)(\alpha_2-\beta_1)[\overline{\alpha}_1^2+2\overline{\alpha}_1(\overline{\alpha}_2-\overline{\alpha}_1)+(\overline{\alpha}_2-\overline{\alpha}_1)^2]+c(\alpha_2-\alpha_1)^2=0,
    \end{split}
    \end{equation*}
which yields equation \eqref{PPreq11srt}. On the other hand, by exchanging the roles of $\alpha_1$ and $\alpha_2$, we also have
    \begin{equation}\label{PPreq11srt ex}
    4\overline{\alpha}_1^{2}(\alpha_1-\beta_1)(\alpha_2-\beta_1)+c(\alpha_2-\alpha_1)^{2}=0.
    \end{equation}
From \eqref{PPreq11srt} and \eqref{PPreq11srt ex}, we see that $\alpha_2=-\alpha_1$ and it follows that $c=\overline{\alpha}_1^2(1-\beta_1^2/\alpha_1^2)$.
Denote $\delta=\beta_1/\alpha_1$. Then $\delta\not=\pm1$ is an algebroid function of degree at most~2.
By doing the transformation $f\rightarrow\alpha_1f$, then \eqref{PPreq8s} becomes
    \begin{equation*}\label{PPreq12s}
    \overline{f}^2=\frac{(1-\delta^2)(f^2-1)}{(f-\delta)^2}=1-\left(\frac{\delta f-1}{f-\delta}\right)^2,
    \end{equation*}
which is the equation \eqref{yanagiharaeq12}.

If $c_2=1/(\alpha_2-\alpha_1)$ is not a root of \eqref{PPreq11s}, then roots of \eqref{PPreq11s} must coincide with $v_1,v_2$, which means that
    \begin{equation}\label{PPreq14sstt}
    \begin{split}
    \frac{\overline{\overline{\alpha}}_2^{2}(\overline{\alpha}_1-\overline{\beta}_1)^2}{4(\alpha_1-\beta_1)(\alpha_2-\beta_1)\overline{\alpha}_1^2+c(\alpha_2-\alpha_1)^2}=\frac{2\overline{\overline{\alpha}}_2^{2}(\overline{\alpha}_1-\overline{\beta}_1)+\overline{c}_1^{2}}{8(\alpha_1-\beta_1)(\alpha_2-\beta_1)\overline{\alpha}_1}=\frac{\overline{\overline{\alpha}}_2^{2}-\overline{c}_1^{2}\overline{c}_2}{4(\alpha_1-\beta_1)(\alpha_2-\beta_1)}.
    \end{split}
    \end{equation}
Recall that $c_1=[c(\alpha_2-\alpha_1)]^{1/2}$ and $c_2=1/(\alpha_2-\alpha_1)$. We obtain from the second equation of \eqref{PPreq14sstt} that
    \begin{equation*}
    (\overline{\alpha}_1+\overline{\alpha}_2)\overline{c}=2\overline{\overline{\alpha}}_2^2\overline{\beta}_1.
     \end{equation*}
By exchanging the roles of $\alpha_1$ and $\alpha_2$, we see that $\overline{\beta}_1\equiv0$ and then it follows from the first equation of \eqref{PPreq14sstt} that
     \begin{equation*}\label{PPreq14ss be}
    \overline{c}c=\overline{c}\overline{\alpha}_1^2+c\overline{\overline{\alpha}}_1^2.
     \end{equation*}
Denote $d=c/\overline{\alpha}_1^2$. Then $d$ is a meromorphic function
satisfying $\overline{d}d=\overline{d}+d$. By doing the transformation $f\rightarrow\alpha_1f$, then \eqref{PPreq8s} becomes
    \begin{equation*}\label{PPreq14ss}
    \overline{f}^2=d(1-f^{-2}),
    \end{equation*}
which is the equation \eqref{yanagiharaeq12a}.

We now look at equation \eqref{Preq3s}. By using a similar analysis on the $\alpha_i$-points of $f$ as to the $\alpha_i$-points of $f$ in \eqref{Preq3} in Section~\ref{proof2_sec}, we obtain the following three possibilities:
    \begin{eqnarray}
    \overline{f}^3 &=& \frac{c(f-\alpha_1)(f-\alpha_2)(f-\alpha_3)}{(f-\beta_1)^3}, \quad \alpha_1\alpha_2\alpha_3\not\equiv0,\label{PQreq24s}\\
    \overline{f}^4 &=& \frac{c(f-\alpha_1)^{2}(f-\alpha_2)(f-\alpha_3)}{(f-\beta_1)^4}, \quad \alpha_1\alpha_2\alpha_3\not\equiv0,\label{PQreq25s}\\
    \overline{f}^6 &=& \frac{c(f-\alpha_1)^{3}(f-\alpha_2)^{2}(f-\alpha_3)}{(f-\beta_1)^6}, \quad \alpha_1\alpha_2\alpha_3\not\equiv0.\label{PQreq26s}
   \end{eqnarray}

Let's first consider equation \eqref{PQreq24s}. Put
    \begin{equation}\label{PQreq27sbefore}
    \frac{\overline{f}(f-\beta_1)}{f-\alpha_1}=c_1u, \quad \frac{1}{f-\alpha_1}=b_2v+b_1,
   \end{equation}
where $u$ and $v$ are two algebroid functions with at most $S(r,f)$ many branch points and the coefficients $c_1,b_1,b_2$ take the same form as in \eqref{coefficients1} which are algebroid functions of degrees at most~$9,3,3$, respectively.
By the First Main Theorem, we have $T(r,v)=T(r,f)+S(r,f)$. Thus all the coefficients of \eqref{PQreq24s} are small functions of $v$. Then we have
    \begin{equation*}
    f=\alpha_1+\frac{1}{b_2v+b_1}, \quad \overline{f}=\frac{c_1u}{1+(\alpha_1-\beta_1)(b_2v+b_1)}
   \end{equation*}
and from \eqref{PQreq24s} and \eqref{PQreq27sbefore} that
    \begin{equation}\label{PQreq27s}
    v^2=4u^3+1.
    \end{equation}
It also follows that
    \begin{equation}\label{PQreq28s}
   \frac{c_1u-\overline{\alpha}_2[1+(\alpha_1-\beta_1)(b_2v+b_1)]}{1+(\alpha_1-\beta_1)(b_2v+b_1)}=\overline{f}-\overline{\alpha}_2.
   \end{equation}
Since zeros of $\overline{f}-\overline{\alpha}_2$ are threefold ramified, we see that zeros of $F=c_1u-\overline{\alpha}_2[1+(\alpha_1-\beta_1)(b_2v+b_1)]$
are also threefold ramified. Let $u_0$ and $v_0$ be a pair of roots of the following equations
    \begin{equation*}
    v_0^2=4u_0^3+1, \quad c_1u_0-\overline{\alpha}_2[1+(\alpha_1-\beta_1)(b_2v_0+b_1)]=0.
    \end{equation*}
Eliminating $u_0$ from the above two equations gives
    \begin{equation}\label{PPreq29yys}
    4\overline{\alpha}_2^3[1+(\alpha_1-\beta_1)(b_2v_0+b_1)]^3-c_1^3(v_0^2-1)=0,
    \end{equation}
from which we see that $v_0$ is an algebroid function of degree at most~$9$.
Moreover, it is seen that~$\pm1$ are not roots of \eqref{PPreq29yys} for otherwise we have $1+(\alpha_1-\beta_1)(b_2v_0+b_1)=0$ and it follows that $\pm b_2+b_1=-1/(\alpha_1-\beta_1)$, which gives $\alpha_2=\beta_1$ or $\alpha_3=\beta_1$, a contradiction.
Rewrite \eqref{PQreq28s} as
   \begin{equation}\label{PQreq28afteryys}
   -\frac{4\overline{\alpha}_2^3[1+(\alpha_1-\beta_1)(b_2v+b_1)]^3-c_1^3(v^2-1)}{4[1+(\alpha_1-\beta_1)(b_2v+b_1)]^3}=\overline{f}^3-\overline{\alpha}_2^3.
   \end{equation}
From equation \eqref{PQreq27s}, we see that $v$ is threefold ramified over~$\pm1,\infty$.
This implies that \eqref{PPreq29yys} can admit only one root for otherwise by a similar analysis as to \eqref{Preq23hh} in Section~\ref{proof2_sec} we obtain from \eqref{PQreq28afteryys} that \eqref{PPreq29yys} has at least two roots, say, $v_1$ and $v_2$, distinct from~$\pm1,\infty$ such that $v$ is completely ramified at $v_1$-points and $v_2$-points, which is impossible by Theorem~\ref{completelyrm}. Expanding \eqref{PPreq29yys} gives
    \begin{equation}\label{PPreq29extrasyys}
        \begin{split}
        &4\overline{\alpha}_2^3(\alpha_1-\beta_1)^3b_2^3v_0^3+12\overline{\alpha}_2^3(\alpha_1-\beta_1)^2b_2^2[b_1(\alpha_1-\beta_1)+1]v_0^2-c_1^3v_0^2\\
        &+12\overline{\alpha}_2^3(\alpha_1-\beta_1)b_2[b_1(\alpha_1-\beta_1)+1]^2v_0+4\overline{\alpha}_2^3[b_1(\alpha_1-\beta_1)+1]^3+c_1^3=0,
        \end{split}
    \end{equation}
and so by Vieta's formulas, we get
    \begin{equation}\label{PPreq29extrasyyseza}
    \frac{12\overline{\alpha}_2^3(\alpha_1-\beta_1)^2b_2^2[b_1(\alpha_1-\beta_1)+1]-c_1^3}{4\overline{\alpha}_2^3(\alpha_1-\beta_1)^3b_2^3}=-3v_0
    \end{equation}
and
     \begin{equation}\label{PPreq29extrasyysezb}
    \frac{3[b_1(\alpha_1-\beta_1)+1]^2}{(\alpha_1-\beta_1)^2b_2^2}=3v_0^2.
    \end{equation}
From \eqref{PPreq29extrasyysezb}, we have
$v_0=\pm[b_1(\alpha_1-\beta_1)+1]/[(\alpha_1-\beta_1)b_2]$. If $v_0=-[b_1(\alpha_1-\beta_1)+1]/[(\alpha_1-\beta_1)b_2]$, then substitution into \eqref{PPreq29yys} yields $v_0^2=1$, a contradiction. Hence $v_0=[b_1(\alpha_1-\beta_1)+1]/[(\alpha_1-\beta_1)b_2]$. Combining this equation with \eqref{PPreq29extrasyyseza}, we get
    \begin{equation*}
    c_1^3=24\overline{\alpha}_2^3(\alpha_1-\beta_1)^2b_2^2[b_1(\alpha_1-\beta_1)+1].
    \end{equation*}
Similarly, the above equation also holds when replacing $\overline{\alpha}_2$ in \eqref{PQreq28s} with $\overline{\alpha}_1$ and $\overline{\alpha}_3$, which implies that $\alpha_1^3=\alpha_2^3=\alpha_3^3$. Thus, $\alpha_2=\eta\alpha_1$, $\alpha_3=\eta^2\alpha_1$, where $\eta$ is the cubic root of~1 satisfying $\eta^2+\eta+1=0$. Moreover, we have
     \begin{equation*}
    c_1^3=4cb_2^2(\alpha_2-\alpha_1)(\alpha_3-\alpha_1)=24\overline{\alpha}_1^3(\alpha_1-\beta_1)^2b_2^2[b_1(\alpha_1-\beta_1)+1].
    \end{equation*}
Denote $\gamma_1=\beta_1/\alpha_1$, $\gamma_1^3\not\equiv1$. Using equation $\eta^2+\eta+1=0$, we get from the above equation that
      \begin{equation}\label{PQreq38sbee1}
    c=\overline{\alpha}_1^3(1-\gamma_1)^2(1+\gamma_1).
    \end{equation}
Denoting $\gamma_2=\beta_1/\alpha_2=\beta_1/(\eta\alpha_1)$, $\gamma_3=\beta_1/\alpha_3=\beta_1/(\eta^2\alpha_1)$, we have $\gamma_2\not\equiv1$, $\gamma_3\not\equiv1$. Exchanging the role of $\alpha_1$ to $\alpha_2$ and $\alpha_3$, respectively, we also have
    \begin{equation*}
    c=\overline{\alpha}_3^3(1-\eta\gamma_1)^2(1+\eta\gamma_1)=\overline{\alpha}_2^3(1-\eta^2\gamma_1)^2(1+\eta^2\gamma_1),
    \end{equation*}
which together with the equation $\eta^2+\eta+1=0$ yields
    \begin{equation*}
   (\eta^2-\eta)\gamma_1(\gamma_1-1)=0.
    \end{equation*}
Therefore, we have $\gamma_1\equiv0$ and so $\beta_1\equiv0$ and then from \eqref{PQreq38sbee1} that $c=\overline{\alpha}_1^3$. By doing the transformation $f\rightarrow\alpha_1f$, then \eqref{PQreq24s} becomes
    \begin{equation*}\label{PQreq38s}
     \overline{f}^3=1-f^{-3},
    \end{equation*}
which is the equation \eqref{yanagiharaeq15}.

We now consider equation \eqref{PQreq25s}. Put
    \begin{equation}\label{PQreq39sbess}
    \frac{\overline{f}(f-\beta_1)}{f-\alpha_1}=c_1u, \quad \frac{1}{f-\alpha_1}=b_2v+b_1,
    \end{equation}
where $u$ and $v$ are two algebroid functions with at most $S(r,f)$ many branch points and the coefficients $c_1,b_1,b_2$ take the same form as in \eqref{coefficients2} which
are algebroid functions of degrees at most~$16,4,4$, respectively.
By the First Main Theorem, we have $T(r,v)=T(r,f)+S(r,f)$. Thus all the coefficients of \eqref{PQreq25s} are small functions of $v$. Then we have
    \begin{equation*}
    f=\frac{1}{b_2v+b_1}+\alpha_1, \quad \overline{f}=\frac{c_1u}{1+(\alpha_1-\beta_1)(b_2v+b_1)}
   \end{equation*}
and from \eqref{PQreq25s} and \eqref{PQreq39sbess} that
    \begin{equation}\label{PQreq39s}
    v^2=1-u^4.
   \end{equation}
It also follows that
    \begin{equation}\label{PQreq39swqes}
   \frac{c_1u-\overline{\alpha}_2[1+(\alpha_1-\beta_1)(b_2v+b_1)]}{1+(\alpha_1-\beta_1)(b_2v+b_1)}=\overline{f}-\overline{\alpha}_2.
   \end{equation}
Since zeros of $\overline{f}-\overline{\alpha}_2$ are fourfold ramified, we see that zeros of $F=c_1u-\overline{\alpha}_2[1+(\alpha_1-\beta_1)(b_2v+b_1)]$ are also fourfold ramified. Let $u_0$ and $v_0$ be a pair of roots of the following equations
    \begin{equation*}
    v_0^2=1-u_0^4, \quad c_1u_0-\overline{\alpha}_2[1+(\alpha_1-\beta_1)(b_2v_0+b_1)]=0.
    \end{equation*}
Eliminating $u_0$ from the above two equations gives
    \begin{equation}\label{PQreq41bbb}
    \overline{\alpha}_2^4[1+(\alpha_1-\beta_1)(b_2v_0+b_1)]^4-c_1^4(1-v_0^2)=0,
    \end{equation}
from which we see that $v_0$ is an algebroid function of degree at most~$16$.
Moreover, it is seen that $\pm1$ are not roots of \eqref{PQreq41bbb} for otherwise we have $1+(\alpha_1-\beta_1)(b_2v_0+b_1)=0$ and it follows that $\pm b_2+b_1=-1/(\alpha_1-\beta_1)$, which gives $\alpha_2=\beta_1$ or $\alpha_3=\beta_1$, a contradiction. Rewrite \eqref{PQreq39swqes} as
    \begin{equation}\label{PQreq39swqess}
   \frac{c_1^4(1-v^2)-\overline{\alpha}_2^4[1+(\alpha_1-\beta_1)(b_2v+b_1)]^4}{[1+(\alpha_1-\beta_1)(b_2v+b_1)]^4}=\overline{f}^4-\overline{\alpha}_2^4.
   \end{equation}
From equation \eqref{PQreq39s}, we see that $v$ is completely ramified over~$\pm1,\infty$.
This implies that \eqref{PQreq41bbb} can admit only one root for otherwise by a similar analysis as to \eqref{Preq23hh} in Section~\ref{proof2_sec} we obtain from \eqref{PQreq39swqess} that
\eqref{PQreq41bbb} has at least two roots, say, $v_1$ and $v_2$, distinct from~$\pm1,\infty$ such that $v$ is completely ramified at $v_1$-points and $v_2$-points, which is impossible by Theorem~\ref{completelyrm}. Expanding \eqref{PQreq41bbb} gives
    \begin{equation*}
    \begin{split}
    \overline{\alpha}_2^4\{(\alpha_1-\beta_1)^4b_2^4v_0^4+4(\alpha_1-\beta_1)^3b_2^3[1+(\alpha_1-\beta_1)b_1]v_0^3+\ldots+[1+(\alpha_1-\beta_1)b_1]^4\}-c_1^4(1-v_0^2)=0,
    \end{split}
    \end{equation*}
and so by Vieta's formulas, we obtain $4[1+(\alpha_1-\beta_1)b_1]/(\alpha_1-\beta_1)b_2=-4v_0$,
which gives $1+(\alpha_1-\beta_1)(b_2v_0+b_1)=0$ and it follows from \eqref{PQreq41bbb} that $v_0^2=1$, a contradiction. Therefore, \eqref{PQreq25s} cannot admit any meromorphic solutions.

Finally, let's consider equation \eqref{PQreq26s}. Put
    \begin{equation}\label{PQreq51sbes}
    \frac{\overline{f}(f-\beta_1)}{f-\alpha_1}=c_1u, \quad
    \frac{1}{f-\alpha_1}=b_2v-b_1,
    \end{equation}
where $u$ and $v$ are two algebroid functions with at most $S(r,f)$ many branch points and the coefficients $c_1,b_1,b_2$ take the same form as in \eqref{coefficients1} which are algebroid functions of degrees at most~$36,6,6$, respectively.
By the First Main Theorem, we have $T(r,v)=T(r,f)+S(r,f)$. Thus all the coefficients of \eqref{PQreq26s} are small functions of $v$. Then we have
    \begin{equation*}
    f=\frac{1}{b_2v-b_1}+\alpha_1, \quad \overline{f}=\frac{c_1u}{1+(\alpha_1-\beta_1)(b_2v-b_1)}
   \end{equation*}
and from \eqref{PQreq26s} and \eqref{PQreq51sbes} that
    \begin{equation}\label{PQreq51s}
    u^6=(v-1)^2(v+2)=v^3-3v+2.
   \end{equation}
It also follows that
    \begin{equation}\label{PQreq51sure}
    \frac{c_1u-\overline{\alpha}_3[1+(\alpha_1-\beta_1)(b_2v-b_1)]}{1+(\alpha_1-\beta_1)(b_2v-b_1)}=\overline{f}-\overline{\alpha}_3.
   \end{equation}
Since zeros of $\overline{f}-\overline{\alpha}_3$ are sixfold ramified, we see that zeros of $F=c_1u-\overline{\alpha}_3[1+(\alpha_1-\beta_1)(b_2v-b_1)]$
are also sixfold ramified. Let $u_0$ and $v_0$ be a pair of roots of the following equations
    \begin{equation*}
    u_0^6=v_0^3-3v_0+2, \quad c_1u_0-\overline{\alpha}_3[1+(\alpha_1-\beta_1)(b_2v_0-b_1)]=0.
    \end{equation*}
Eliminating $u_0$ from the above two equations gives
    \begin{equation}\label{PQreq53s}
    \overline{\alpha}_3^6[1+(\alpha_1-\beta_1)(b_2v_0-b_1)]^6-c_1^6(v_0^3-3v_0+2)=0,
    \end{equation}
from which we see that $v_0$ is an algebroid function of degree at most~$36$.
Moreover, it is easy to see that~$1$ and~$-2$ are not roots of \eqref{PQreq53s} for otherwise we have $1+(\alpha_1-\beta_1)(b_2-b_1)=0$, which gives $\alpha_2=\beta_1$, a contradiction; or $1+(\alpha_1-\beta_1)(-2b_2-b_1)=0$, which gives $\alpha_3=\beta_1$, a contradiction.
Rewrite \eqref{PQreq51sure} as
    \begin{equation}\label{PQreq51sure'}
    \frac{c_1^6(v^3-3v+2)-\overline{\alpha}_3^6[1+(\alpha_1-\beta_1)(b_2v-b_1)]^6}{[1+(\alpha_1-\beta_1)(b_2v-b_1)]^6}=\overline{f}^6-\overline{\alpha}_3^6.
   \end{equation}
From equation \eqref{PQreq51s}, we see that $v$ is completely ramified over~$1,-2,\infty$.
This implies that \eqref{PQreq53s} can admit only one root for otherwise by a similar analysis as to \eqref{Preq23hh} in Section~\ref{proof2_sec} we obtain from \eqref{PQreq51sure'} that
\eqref{PQreq53s} has at least two roots, say, $v_1$ and $v_2$, distinct from~$1,-2,\infty$ such that $v$ is completely ramified at $v_1$-points and $v_2$-points, which is impossible by Theorem~\ref{completelyrm}. Expanding \eqref{PQreq53s} gives
     \begin{equation*}
     \begin{split}
     \overline{\alpha}_3^6\{(\alpha_1-\beta_1)^6b_2^6v_0^6+6(\alpha_1-\beta_1)^5b_2^5v_0^5[1-(\alpha_1-\beta_1)b_1]+\ldots+[1-(\alpha_1-\beta_1)b_1]^6\}-c_1^6(v_0^3-3v_0+2)=0,
    \end{split}
    \end{equation*}
and so by Vieta's formulas, we obtain $-6[1-(\alpha_1-\beta_1)b_1]/(\alpha_1-\beta_1)b_2=-6v_0$,
which gives $1+(\alpha_1-\beta_1)(b_2v_0-b_1)=0$ and it follows from \eqref{PQreq53s} that $v_0^3-3v_0+2=0$, a contradiction. Therefore, \eqref{PQreq26s} cannot admit any meromorphic solutions.

\subsection{$Q(z,f)$ has two distinct roots}\label{proof32_sec}

Now $P(z,f)$ can have one or two roots. We discuss these two cases separately as follows:

\vskip 0.3cm

\noindent{\bf Case 1:} $P(z,f)$ has only one root. Then \eqref{yanagiharaeq2} assumes the following form:
    \begin{equation}\label{SSS1yy}
    \overline{f}^{n}=\frac{c(f-\alpha_1)^{n}}{(f-\beta_1)^{\nu_1}(f-\beta_2)^{\nu_2}},
    \end{equation}
where $\nu_1,\nu_2$ are positive integers satisfying $\nu_1+\nu_2=n$. By doing a bilinear transformation $f\rightarrow 1/f$, if $\alpha_1\equiv0$, we obtain \eqref{Preq2}; if $\alpha_1\not\equiv0$ and $\beta_1,\beta_2\not\equiv0$, we obtain \eqref{Preq2s}.
So in what follows we only consider the case that $\alpha_1\not\equiv0$ and $\beta_1\equiv0$.
Note that $(\nu_1,\nu_2)=1$. If $n=2$, then obviously $\nu_1=\nu_2=1$; if $n\geq3$, then by using a similar analysis on the $\beta_j$-points of $f$ in \eqref{SSS1yy} as to $\alpha_i$-points of $f$ in \eqref{Preq2} in Section~\ref{proof2_sec}, we can obtain that $\nu_1=1$ or $\nu_2=1$. Thus \eqref{SSS1yy} is one of the following forms:
     \begin{eqnarray}
    \overline{f}^{n}&=&\frac{c(f-\alpha_1)^{n}}{f(f-\beta_2)^{n-1}},\label{SSS1}\\
    \overline{f}^{n}&=&\frac{c(f-\alpha_1)^{n}}{f^{n-1}(f-\beta_2)}.\label{SSS188}
    \end{eqnarray}
Consider equation \eqref{SSS1} first. By doing a bilinear transformation $f\rightarrow1/f$, we obtain
     \begin{equation}\label{SSS2}
    \overline{f}^{n}=\frac{d_1(f-\gamma_1)^{n-1}}{(f-\lambda_1)^{n}},
    \end{equation}
where $d_1=-\beta_2^{n-1}/(c\alpha_1^n)$, $\gamma_1=1/\beta_2$ and $\lambda_1=1/\alpha_1$ are nonzero small meromorphic functions of $f$.
Moreover, it is seen that $f$ is completely ramified at poles, $\gamma_1$-points and also $\lambda_1$-points. Denote $c_1=d_1^{1/n}$. Then $c_1$ is an algebroid function of degree at most~$n$. Put
    \begin{equation}\label{Preq20skk}
    c_1u=\frac{\overline{f}(f-\lambda_1)}{f-\gamma_1}, \quad v=\frac{1}{f-\gamma_1},
    \end{equation}
where $u$ is an algebroid function with at most $S(r,f)$ many branch points and $v$ is a meromorphic function.
By the First Main Theorem, we have $T(r,v)=T(r,f)+S(r,f)$. Thus all the coefficients of \eqref{SSS1} are small functions of $v$. Then we have
    \begin{equation*}
    f=\frac{1}{v}+\gamma_1, \quad \overline{f}=\frac{c_1u}{1+(\gamma_1-\lambda_1)v}
    \end{equation*}
and from \eqref{SSS2} and \eqref{Preq20skk} that
    \begin{equation}\label{Preq21skk}
    u^{n}=v.
    \end{equation}
It also follows that
    \begin{equation}\label{Preq21skk'}
    \frac{c_1u-\overline{\gamma}_1[1+(\gamma_1-\lambda_1)v]}{1+(\gamma_1-\lambda_1)v}=\overline{f}-\overline{\gamma}_1.
    \end{equation}
Since $\overline{f}$ is completely ramified at $\overline{\gamma}_1$-points, we see that $F=c_1u-\overline{\gamma}_1[(\gamma_1-\lambda_1)v+1]$ is completely ramified at zeros. Let $u_0$ and $v_0$ be a pair of roots of the following equations
    \begin{equation*}
    u_0^{n}=v_0, \quad c_1u_0-\overline{\gamma}_1[(\gamma_1-\lambda_1)v_0+1]=0.
    \end{equation*}
Eliminating $u_0$ from the above two equations gives
    \begin{equation*}
    \overline{\gamma}_1^{n}[(\gamma_1-\lambda_1)v_0+1]^{n}-c_1^{n}v_0=0,
    \end{equation*}
i.e.,
    \begin{equation}\label{PPreq6skk}
    \overline{\gamma}_1^{n}[(\gamma_1-\lambda_1)^nv_0^n+n(\gamma_1-\lambda_1)^{n-1}v_0^{n-1}+\ldots+1]-c_1^{n}v_0=0,
    \end{equation}
from which we see that $v_0$ is an algebroid function of degree at most $n$. Moreover, it is seen that~0 and $1/(\lambda_1-\gamma_1)$ are not roots of \eqref{PPreq6skk}. Rewrite \eqref{Preq21skk'} as
    \begin{equation}\label{Preq21skk''}
    \frac{c_1^nv-\overline{\gamma}_1^n[1+(\gamma_1-\lambda_1)v]^n}{[1+(\gamma_1-\lambda_1)v]^n}=\overline{f}^n-\overline{\gamma}_1^n.
    \end{equation}
We see from \eqref{Preq20skk} that $v$ is completely ramified at zeros, poles and $1/(\lambda_1-\gamma_1)$-points. This implies that \eqref{PPreq6skk} can admit only one root for otherwise by using a similar analysis as to \eqref{Preq23hh} in Section~\ref{proof2_sec} we obtain from \eqref{Preq21skk''} that \eqref{PPreq6skk} has at least two roots, say, $v_1$ and $v_2$, distinct from~$0,\infty,1/(\lambda_1-\gamma_1)$ such that $v$ is completely ramified at $v_1$-points and $v_2$-points, which is impossible by Theorem~\ref{completelyrm}. However, when $n\geq 3$, by Vieta's formulas we see that the only root is $v_0=-1/(\gamma_1-\lambda_1)$ and it follows from \eqref{PPreq6skk} that $c_1^{n}/(\gamma_1-\lambda_1)=0$, a contradiction. Hence $n=2$ and we therefore have
     \begin{equation}\label{SSS2des}
    \overline{f}^{2}=\frac{c(f-\gamma_1)}{(f-\lambda_1)^{2}}.
    \end{equation}
Moreover, we conclude from the above reasoning that the quadratic equation
    \begin{equation}\label{PPreq6skkdes}
    \overline{\gamma}_1^{2}[(\gamma_1-\lambda_1)v_0+1]^{2}-c_1^{2}v_0=0
    \end{equation}
can admit only one root. It follows that the discriminant of equation \eqref{PPreq6skkdes} with respect to $v_0$ satisfies
    \begin{equation*}
    [2\overline{\gamma}_1^{2}(\gamma_1-\lambda_1)-c_1^{2}]^2-4\overline{\gamma}_1^{4}(\gamma_1-\lambda_1)^2=0,
    \end{equation*}
i.e.,
    \begin{equation}\label{PPreq8skk}
    d_1=4\overline{\gamma}_1^{2}(\gamma_1-\lambda_1).
    \end{equation}
On the other hand, since $\overline{f}$ is completely ramified at $\overline{\lambda}_1$-points, if we consider equation
    \begin{equation*}
    \frac{c_1u-\overline{\lambda}_1[1+(\gamma_1-\lambda_1)v]}{1+(\gamma_1-\lambda_1)v}=\overline{f}-\overline{\lambda}_1,
    \end{equation*}
then by similar arguments as above we also have
    \begin{equation}\label{PPreq9skk}
    d_1=4\overline{\lambda}_1^{2}(\gamma_1-\lambda_1).
    \end{equation}
From \eqref{PPreq8skk} and \eqref{PPreq9skk}, we see that $\lambda_1=-\gamma_1$ and it follows that $d_1=-8\overline{\lambda}_1^{2}\lambda_1$. By doing the transformation $f\rightarrow\gamma_1f$, then \eqref{SSS2des} becomes
     \begin{equation*}
     \overline{f}^{2}=-\frac{8(f+1)}{(f-1)^{2}}=1-\left(\frac{f+3}{f-1}\right)^{2},
    \end{equation*}
which is the equation \eqref{yanagiharaeq13}.

Consider now equation \eqref{SSS188}. By doing a bilinear transformation $f\rightarrow1/f$, then \eqref{SSS188} becomes
     \begin{equation}\label{SSS2re}
    \overline{f}^{n}=\frac{d_2(f-\gamma_2)}{(f-\lambda_2)^{n}},
    \end{equation}
where $d_2=(-1)^{1-n}\beta_1/(c\alpha_1^n)$, $\gamma_2=1/\beta_1$ and $\lambda_2=1/\alpha_1$ are nonzero small meromorphic functions of $f$. Moreover, $f$ is completely ramified at poles, $\gamma_2$-points and also $\lambda_2$-points. By a similar reasoning as above, we can also obtain $n=2$ and this will lead \eqref{SSS2re} into \eqref{yanagiharaeq13} again. We omit all those details.

\vskip 0.3cm

\noindent{\bf Case 2:} $P(z,f)$ has two distinct roots. Then \eqref{yanagiharaeq2} assumes the following form:
    \begin{equation}\label{SSSD}
    \overline{f}^{n}=\frac{c(f-\alpha_1)^{\kappa_1}(f-\alpha_2)^{\kappa_2}}{(f-\beta_1)^{\nu_1}(f-\beta_2)^{\nu_2}}, \quad  \alpha_1\alpha_2\not\equiv0,
    \end{equation}
where $\kappa_i$, $\nu_j$ are positive integers satisfying $\kappa_1+\kappa_2=n$, $\nu_1+\nu_2=n$.
Denote $k_1=(n,\kappa_1)$, $k_2=(n,\kappa_2)$, $l_1=(n,\nu_1)$, $l_2=(n,\nu_2)$ and let $K_1=n/k_1$, $K_2=n/k_2$, $L_1=n/l_1$, $L_2=n/l_2$. Obviously, $K_1\geq2$, $K_2\geq2$, $L_1\geq2$, $L_2\geq2$.
Now $f$ is $K_i$-fold ramified at $\alpha_i$-points and $L_j$-fold ramified at $\beta_j$-points, by \eqref{de relation}, we conclude that $K_1=K_2=L_1=L_2=2$ and $k_1+k_2=n$, $l_1+l_2=n$. Since $\kappa_1+\kappa_2=n$, $\nu_1+\nu_2=n$, we see that $\kappa_1=\kappa_2=\nu_1=\nu_2=n/2$. Hence $n=2$ and we therefore have
    \begin{equation}\label{SSSD1}
    \overline{f}^2=\frac{c(f-\alpha_1)(f-\alpha_2)}{(f-\beta_1)(f-\beta_2)},  \quad  \alpha_1\alpha_2\not\equiv0.
    \end{equation}
Moreover, if $\beta_1\equiv0$ or $\beta_2\equiv0$, then $f$ will be fourfold ramified at $\alpha_i$-points, which yields a contradiction to \eqref{de relation}. Therefore, $\beta_1\not\equiv0$ and $\beta_2\not\equiv0$.

Denote $c_1=c^{1/2}$. Then $c_1$ is an algebroid function of degree at most~2. Put
    \begin{equation}\label{SSSD2before}
    c_1u=\frac{\overline{f}(f-\beta_1)}{f-\alpha_1}, \quad v=\frac{f-\beta_1}{f-\alpha_1},
    \end{equation}
where $u$  and $v$ are two algebroid functions with at most $S(r,f)$ many branch points. By the First Main Theorem, we have $T(r,v)=T(r,f)+S(r,f)$. Thus all the coefficients of \eqref{SSSD1} are small functions of $v$. Then we have
    \begin{equation*}
    f=\frac{\alpha_1v-\beta_1}{v-1}, \quad \overline{f}=\frac{c_1u}{v}
    \end{equation*}
and from \eqref{SSSD1} and \eqref{SSSD2before} that
    \begin{equation}\label{SSSD2}
    u^2=v\frac{(\alpha_1-\alpha_2)v-(\beta_1-\alpha_2)}{(\alpha_1-\beta_2)v-(\beta_1-\beta_2)}.
    \end{equation}
It also follows that
   \begin{equation}\label{SSSD3}
    \frac{c_1u-\overline{\alpha}_1v}{v}=\overline{f}-\overline{\alpha}_1.
    \end{equation}
Since zeros of $\overline{f}-\overline{\alpha}_1$ are twofold ramified, we see that zeros of $F=c_1u-\overline{\alpha}_1v$ are also twofold ramified. Let $u_0$ and $v_0$ be a pair of roots of the following equations
    \begin{equation*}
    u_0^2=v_0\frac{(\alpha_1-\alpha_2)v_0-(\beta_1-\alpha_2)}{(\alpha_1-\beta_2)v_0-(\beta_1-\beta_2)}, \quad c_1u_0-\overline{\alpha}_1v_0=0.
    \end{equation*}
Eliminating $u_0$ from the above two equations gives $v_0\equiv 0$ or
    \begin{equation}\label{SSSD4}
    \overline{\alpha}_1^2(\alpha_1-\beta_2)v_0^2+c_1^2(\beta_1-\alpha_2)-[\overline{\alpha}_1^2(\beta_1-\beta_2)+c_1^2(\alpha_1-\alpha_2)]v_0=0,
    \end{equation}
from which we see that $v_0$ is an algebroid function of degree at most $4$. Moreover, it is seen that none of~0, $(\alpha_2-\beta_1)/(\alpha_2-\alpha_1)$ and $(\beta_2-\beta_1)/(\beta_2-\alpha_1)$ solves \eqref{SSSD4} for otherwise we get $\overline{\alpha}_1\equiv0$, a contradiction. If $v_0\equiv0$, then for a point $z_0$ such that $f(z_0+1)=\alpha_1(z_0+1)$, we have $u(z_0)=v(z_0)=0$ and from \eqref{SSSD2} we see that $c_1(z_0)u(z_0)/v(z_0)=\infty$, a contradiction since $\overline{f}=c_1u/v$. Therefore, $v_0\not\equiv0$.
Rewrite \eqref{SSSD3} as
   \begin{equation}\label{SSSD3''}
    \frac{c_1^2[(\alpha_1-\alpha_2)v-(\beta_1-\alpha_2)]v-[(\alpha_1-\beta_2)v-(\beta_1-\beta_2)]\overline{\alpha}_1^2v^2}{[(\alpha_1-\beta_2)v-(\beta_1-\beta_2)]v^2}=\overline{f}^2-\overline{\alpha}_1^2.
    \end{equation}
From \eqref{SSSD2before} we see that $v$ is completely ramified at poles, zeros, $(\alpha_2-\beta_1)/(\alpha_2-\alpha_1)$-points and $(\beta_2-\beta_1)/(\beta_2-\alpha_1)$-points.
This implies that \eqref{SSSD4} can admit only one root for otherwise by using a similar analysis as to \eqref{Preq23hh} in Section~\ref{proof2_sec} we obtain from \eqref{SSSD3''} that \eqref{SSSD4} has at least two roots, say, $v_1$ and $v_2$, distinct from~$0,\infty,(\alpha_2-\beta_1)/(\alpha_2-\alpha_1),(\beta_2-\beta_1)/(\beta_2-\alpha_1)$ such that $v$ is completely ramified at $v_1$-points and $v_2$-points, which is impossible by Theorem~\ref{completelyrm}.
Therefore, the discriminant of the quadratic equation \eqref{SSSD4} with respect to $v_0$ satisfies
   \begin{equation*}
    [\overline{\alpha}_1^2(\beta_1-\beta_2)+c_1^2(\alpha_1-\alpha_2)]^2-4\overline{\alpha}_1^2c_1^2(\alpha_1-\beta_2)(\beta_1-\alpha_2)=0,
    \end{equation*}
i.e.,
    \begin{equation}\label{SSSD10}
    \overline{\alpha}_1^4-\frac{2c[2(\beta_1-\alpha_2)(\alpha_1-\beta_2)-(\beta_1-\beta_2)(\alpha_1-\alpha_2)]}{(\beta_1-\beta_2)^2}\overline{\alpha}_1^2+\frac{c^{2}(\alpha_1-\alpha_2)^2}{(\beta_1-\beta_2)^2}=0.
    \end{equation}
Similarly, we obtain that \eqref{SSSD10} still holds when replacing $\overline{\alpha}_1$ by any of $\overline{\alpha}_2$, $\overline{\beta}_1$ and $\overline{\beta}_2$, which implies that $\overline{\alpha}_1^2=\overline{\alpha}_2^2$ or $\overline{\alpha}_1^2=\overline{\beta}_1^2$. Thus we need to consider the following two cases: (1), $\alpha_1=-\alpha_2$, $\beta_1=-\beta_2$; and (2), $\alpha_1=-\beta_1$, $\alpha_2=-\beta_2$.

If $\alpha_1=-\alpha_2$, $\beta_1=-\beta_2$, then by Vieta's formulas we obtain from \eqref{SSSD10} that
    \begin{equation*}
    \overline{\alpha}_1^2\overline{\beta}_1^2=\frac{c^{2}\alpha_1^2}{\beta_1^2}, \quad \overline{\alpha}_1^2+\overline{\beta}_1^2=c\frac{\beta_1^2+\alpha_1^2}{\beta_1^2}.
    \end{equation*}
Denote $\gamma_1=\alpha_1/\beta_1$, $\gamma_1^2\not\equiv0,1$ and $\gamma_2=\theta\overline{\gamma}_1/\gamma_1=c/\overline{\beta}_1^2$, $\theta=\pm1$. Then $\gamma_1^2$ is meromorphic and $\gamma_2$ satisfies
$(\overline{\gamma}_1^2+1)\gamma_1=\theta\overline{\gamma}_1(\gamma_1^2+1)$, which gives $\overline{\gamma}_1=\theta\gamma_1$ or $\overline{\gamma}_1\gamma_1=\theta$. It follows that $\gamma_2\equiv1$ or $\gamma_2=\overline{\gamma}_1^2$. If $\gamma_2\equiv1$, then $\gamma_1^2$ is a periodic function with period~1. By doing the transformation $f\rightarrow\beta_1f$, then \eqref{SSSD1} becomes
    \begin{equation}\label{SSSD11}
    \overline{f}^2=\frac{f^2-\gamma_1^2}{f^2-1},
    \end{equation}
which is the equation \eqref{yanagiharaeq13a}; if $\gamma_2=\overline{\gamma}_1^2$, then $\overline{\gamma}_1^2\gamma_1^2=1$. By doing the transformation $f\rightarrow\beta_1f$, then \eqref{SSSD1} becomes
    \begin{equation}\label{SSSD12}
    \overline{f}^2=\frac{\overline{\gamma}_1^2f^2-1}{f^2-1},
    \end{equation}
which is the equation \eqref{yanagiharaeq13b}.

If $\alpha_1=-\beta_1$, $\alpha_2=-\beta_2$, then by Vieta's formulas we obtain from \eqref{SSSD10} that
   \begin{equation*}
    \overline{\alpha}_1^2\overline{\alpha}_2^2=c^2, \quad \overline{\alpha}_1^2+\overline{\alpha}_2^2=2c\frac{-2(\alpha_1+\alpha_2)^2+(\alpha_1-\alpha_2)^2}{(\alpha_1-\alpha_2)^2}.
    \end{equation*}
Denote $\lambda=(\alpha_1/\alpha_2)^{1/2}$, $\lambda^2\not\equiv0,\pm1$. Let $d=\lambda+1/\lambda$. Then $d^2=(\alpha_1+\alpha_2)^{2}/\alpha_1\alpha_2$ is a meromorphic function satisfying
    \begin{equation*}
     \overline{d}^2-2=2\theta\frac{-2d^2+d^2-4}{d^2-4},
    \end{equation*}
i.e.,
    \begin{equation*}
     \overline{d}^2(d^2-4)=2(1-\theta)d^2-8(1+\theta),
    \end{equation*}
where $\theta=\pm1$. By doing the transformation $f\rightarrow(\alpha_1\alpha_2)^{1/2}f$, then \eqref{SSSD1} becomes
    \begin{equation}\label{SSSD14}
    \overline{f}^2=\theta\frac{(f-\lambda)(f-\lambda^{-1})}{(f+\lambda)(f+\lambda^{-1})}=\theta\frac{f^2-df+1}{f^2+df+1},
    \end{equation}
which is the equation \eqref{yanagiharaeq13c}.

\subsection{$Q(z,f)$ has three distinct roots}\label{proof33_sec}

In this case, $f$ is completely ramified at $\beta_j$-points. By using a similar analysis on the $\beta_j$-points of $f$ as to the $\alpha_i$-points of $f$ in \eqref{Preq3} in Section~\ref{proof2_sec}, we obtain the following three possibilities:
    \begin{eqnarray}
    \overline{f}^3 &=& \frac{c(f-\alpha_1)^3}{(f-\beta_1)(f-\beta_2)(f-\beta_3)},\label{TTT1}\\
    \overline{f}^4 &=& \frac{c(f-\alpha_1)^4}{(f-\beta_1)^{2}(f-\beta_2)(f-\beta_3)},\label{TTT2}\\
    \overline{f}^6 &=& \frac{c(f-\alpha_1)^6}{(f-\beta_1)^{3}(f-\beta_2)^{2}(f-\beta_3)},\label{TTT3}
   \end{eqnarray}
where $\alpha_1$ is meromorphic. In \eqref{TTT1}, $f$ is threefold ramified at $\beta_j$-points, $j=1,2,3$. This implies that $\beta_j\not\equiv0$, $j=1,2,3$, for otherwise $f$ will also be threefold ramified at $\alpha_1$-points, which yields a contradiction to \eqref{de relation}. Similarly, in equations \eqref{TTT2}--\eqref{TTT3}, we have $\beta_j\not\equiv0$, $j=1,2,3$. By doing a bilinear transformation $f\rightarrow1/f$, then if $\alpha_1\equiv0$ we obtain equations \eqref{PQreq24}--\eqref{PQreq26} and if $\alpha_1\not\equiv0$ we obtain equations \eqref{PQreq24s}--\eqref{PQreq26s}. This completes the proof of Theorem~\ref{dtheorem}.

%
%
\end{document}